\pdfoutput=1
\documentclass[submission]{FPSAC2024}

\usepackage{tikz}
\usepackage{tkz-euclide} 
\usepackage{ stmaryrd }
\usepackage{ mathrsfs }

\newtheorem{thm}{Theorem}
\newtheorem{lem}[thm]{Lemma}
\newtheorem{prop}[thm]{Proposition}

\theoremstyle{remark}
\newtheorem*{rmk}{Remark}

\newcommand{\GK}{\operatorname{GK}}
\newcommand{\Fer}{\operatorname{Fer}}
\newcommand{\rep}{\operatorname{rep}}
\newcommand{\Mult}{\operatorname{Mult}}
\newcommand{\GenJF}{\operatorname{GenJF}}

\newcommand{\RSK}{\operatorname{RSK}}
\newcommand{\AR}{\operatorname{AR}}
\newcommand{\wt}{\operatorname{\mathsf{wt}}}
\newcommand{\Supp}{\operatorname{Supp}}

\newcommand{\B}{\mathbf B}
\newcommand{\E}{\mathbf E}
\usepackage{lipsum}

\newcommand{\new}[1]{\textit{\color{rub-grun}{#1}}}

\title[RSK via the combinatorics of quiver representations]{An extended generalization of RSK via the combinatorics of type $A$ quiver representations}

\author{Benjamin Dequêne\thanks{dequene.benjamin@courrier.uqam.ca}}

\address{Laboratoire d'Algèbre de Combinatoire et d'Informatique Mathématique (UQAM, Montréal)\\ Laboratoire Amiénois de Mathématiques Fondamentales et Appliquées (UPJV, Amiens)}

\received{\today}

\abstract{The classical Robinson--Schensted--Knuth correspondence is a bijection from nonnegative integer matrices to pairs of semi-standard Young tableaux. Based on the work of, among others, Burge, Hillman, Grassl, Knuth and Gansner, it is known that a version of this correspondence gives, for any nonzero integer partition $\lambda$, a bijection from arbitrary fillings of $\lambda$ to reverse plane partitions of shape $\lambda$, via Greene--Kleitman invariants. By bringing out the combinatorial aspects of our recent results on quiver representations, we construct a family of bijections from fillings of $\lambda$ to reverse plane partitions of shape $\lambda$ parametrized by a choice of Coxeter element in a suitable symmetric group. We recover the above version of the Robinson--Schensted--Knuth correspondence for a particular choice of Coxeter element depending on $\lambda$.}

\resume{La correspondance Robinson--Schensted--Knuth classique est une bijection partant des matrices à coefficients des entiers naturels vers les paires de tableaux de Young semi-standards. Basé sur les travaux, entre autres, de Burge, Hillman, Grassl, Knuth et Gansner, on sait qu'une version de cette correspondance donne, pour toute partage d'un entier non nulle $\lambda$, une bijection allant des remplissages arbitraires de $\lambda$ vers les partitions planes renversées de forme $\lambda$, via les invariants de Greene--Kleitman. En faisant ressortir les aspects combinatoires de nos récents résultats sur les représentations de carquois, nous construisons une famille de bijections partant des remplissages de $\lambda$ vers les partitions planes renversées de forme $\lambda$, paramétrées par un choix d'élément de Coxeter dans un groupe symétrique approprié. Nous récupérons la version de la correspondance Robinson--Schensted--Knuth ci-dessus pour un choix particulier d'élément de Coxeter dépendant de $\lambda$.}


\usepackage[backend=bibtex]{biblatex}
\begin{document}

\maketitle

\section{Introduction} 

The Robinson--Schensted--Knuth (RSK) correspondence is a fundamental bijection from nonnegative integer matrices to pairs of semi-standard Young tableaux of the same shape. For further details, the reader may consult the following references: \cite{St99}, \cite{F96}. 

Based on observations of various works of Burge \cite{B72}, Hillman--Grassl \cite{HG76} and Knuth \cite{K70}, Gansner \cite{Ga81Ma,Ga81Hi} constructed a generalized version of this correspondence, via Greene--Kleitman invariants, which gives a bijection from arbitrary fillings to reverse plane partitions of the same shape. 

Our paper \cite{Deq23} studies a representation-theoretic setting in which a version of RSK exists. In the present paper, we present an explicit, combinatorial form of the results from \cite{Deq23}. Given a fixed nonzero integer partition $\lambda$, we present the construction of a family of maps $(\RSK_{\lambda,c})_c$ from fillings of $\lambda$ to reverse plane partitions of shape $\lambda$ parametrized by $c$ a Coxeter element of the symmetric group $\mathfrak{S}_n$ where $n-1$ is the hook-length of the box $(1,1)$ in $\lambda$. We can state the following result from \cite{Deq23}.

\begin{thm} \label{1stmainthm}
	The map $\RSK_{\lambda,c}$ gives a one-to-one correspondence from fillings of shape $\lambda$ to reverse plane partitions of shape $\lambda$. Moreover, for any $\lambda$, there exists a unique (up to inverse) choice of $c$ such that $\RSK_{\lambda,c}$ coincides with the usual $\RSK$. 
\end{thm}

No knowledge in quiver representation is required to read this abstract, except for \cref{s:quiver} in which we discuss the connection with quiver representations.

\section{Gansner's Ferrers Diagram RSK}
\label{s:Gans}

In this section, we describe Gansner's correspondence explicitly.

\subsection{Some vocabulary}

An \new{integer partition} is a weakly decreasing nonnegative integer sequence $\lambda = (\lambda_n)_{n \in \mathbb{N}^*}$ with finitely many nonzero terms. The \new{length} of $\lambda$ is the minimal $k \in \mathbb{N}$ such that $\lambda_{k+1}=0$. We endow $(\mathbb{N}^*)^2$ with the Cartesian product order $\unlhd$. The \new{Ferrers diagram of $\lambda$} $\Fer(\lambda)$ is the subset of $(\mathbb{N}^*)^2$ given by pairs $(i,j)$ such that $i \leqslant \lambda_j$. We call any map $f : \Fer(\lambda) \longrightarrow \mathbb{N}$ a \new{filling of shape $\lambda$} . Such a filling $f$ is a \new{reverse plane partition} whenever $f$ weakly increases with respect to $\unlhd$. We give an example of a reverse plane partition of shape $(5,3,3,2)$ in \cref{fig:RPP}.
\begin{figure}[h!]
	\centering
	\scalebox{0.5}{\begin{tikzpicture}
			\tkzDefPoint(0,0){a}
			\tkzDefPoint(0,1){b}
			\tkzDefPoint(1,1){c}
			\tkzDefPoint(1,0){d}
			\tkzDrawPolygon[line width = 0.7mm, color = black](a,b,c,d);
			
			\tkzDefPoint(1,0){a}
			\tkzDefPoint(1,1){b}
			\tkzDefPoint(2,1){c}
			\tkzDefPoint(2,0){d}
			\tkzDrawPolygon[line width = 0.7mm, color = black](a,b,c,d);
			
			\tkzDefPoint(2,0){a}
			\tkzDefPoint(2,1){b}
			\tkzDefPoint(3,1){c}
			\tkzDefPoint(3,0){d}
			\tkzDrawPolygon[line width = 0.7mm, color = black](a,b,c,d);
			
			\tkzDefPoint(3,0){a}
			\tkzDefPoint(3,1){b}
			\tkzDefPoint(4,1){c}
			\tkzDefPoint(4,0){d}
			\tkzDrawPolygon[line width = 0.7mm, color = black](a,b,c,d);
			
			\tkzDefPoint(4,1){a}
			\tkzDefPoint(4,0){b}
			\tkzDefPoint(5,0){c}
			\tkzDefPoint(5,1){d}
			\tkzDrawPolygon[line width = 0.7mm, color = black](a,b,c,d);
			
			\tkzDefPoint(0,0){a}
			\tkzDefPoint(0,-1){b}
			\tkzDefPoint(1,-1){c}
			\tkzDefPoint(1,0){d}
			\tkzDrawPolygon[line width = 0.7mm, color = black](a,b,c,d);
			
			\tkzDefPoint(1,0){a}
			\tkzDefPoint(1,-1){b}
			\tkzDefPoint(2,-1){c}
			\tkzDefPoint(2,0){d}
			\tkzDrawPolygon[line width = 0.7mm, color = black](a,b,c,d);
			
			\tkzDefPoint(2,0){a}
			\tkzDefPoint(2,-1){b}
			\tkzDefPoint(3,-1){c}
			\tkzDefPoint(3,0){d}
			\tkzDrawPolygon[line width = 0.7mm, color = black](a,b,c,d);
			
			\tkzDefPoint(0,-2){a}
			\tkzDefPoint(0,-1){b}
			\tkzDefPoint(1,-1){c}
			\tkzDefPoint(1,-2){d}
			\tkzDrawPolygon[line width = 0.7mm, color = black](a,b,c,d);
			
			\tkzDefPoint(1,-2){a}
			\tkzDefPoint(1,-1){b}
			\tkzDefPoint(2,-1){c}
			\tkzDefPoint(2,-2){d}
			\tkzDrawPolygon[line width = 0.7mm, color = black](a,b,c,d);
			
			\tkzDefPoint(2,-2){a}
			\tkzDefPoint(2,-1){b}
			\tkzDefPoint(3,-1){c}
			\tkzDefPoint(3,-2){d}
			\tkzDrawPolygon[line width = 0.7mm, color = black](a,b,c,d);
			
			\tkzDefPoint(0,-3){a}
			\tkzDefPoint(0,-2){b}
			\tkzDefPoint(1,-2){c}
			\tkzDefPoint(1,-3){d}
			\tkzDrawPolygon[line width = 0.7mm, color = black](a,b,c,d);
			
			\tkzDefPoint(1,-3){a}
			\tkzDefPoint(1,-2){b}
			\tkzDefPoint(2,-2){c}
			\tkzDefPoint(2,-3){d}
			\tkzDrawPolygon[line width = 0.7mm, color = black](a,b,c,d);
			
			\tkzLabelPoint[below](0.5,0.9){{\Huge $0$}};
			\tkzLabelPoint[below](1.5,0.9){{\Huge $3$}};
			\tkzLabelPoint[below](2.5,0.9){{\Huge $5$}};
			\tkzLabelPoint[below](3.5,0.9){{\Huge $5$}};
			\tkzLabelPoint[below](4.5,0.9){{\Huge $7$}};
			\tkzLabelPoint[below](0.5,-0.1){{\Huge $1$}};
			\tkzLabelPoint[below](1.5,-0.1){{\Huge $5$}};
			\tkzLabelPoint[below](2.5,-0.1){{\Huge $5$}};
			\tkzLabelPoint[below](0.5,-1.1){{\Huge $4$}};
			\tkzLabelPoint[below](1.5,-1.1){{\Huge $6$}};
			\tkzLabelPoint[below](2.5,-1.1){{\Huge $9$}};
			\tkzLabelPoint[below](0.5,-2.1){{\Huge $4$}};
			\tkzLabelPoint[below](1.5,-2.1){{\Huge $10$}};
	\end{tikzpicture}}
	\caption{\label{fig:RPP} A reverse plane partitions of shape $\lambda = (5,3,3,2)$.}
\end{figure}

\subsection{Greene--Kleitman invariants}

Let $G=(G_0,G_1)$ be a finite directed graph, where $G_0$ is the set of vertices of $G$, and $G_1 \subset (G_0)^2$ is the set of arrows of $G$. Assume that $G$ has no multi-arrows. 

We see a path $\gamma$ in $G$ as a finite sequence of vertices $(v_0, \ldots, v_k)$ such that $(v_i,v_{i+1}) \in G_1$. Denote by $s(\gamma) = v_0$ its source and by $t(\gamma) = v_k$ its target. Write $\Supp(\gamma) = \{v_0, \ldots, v_k\}$ to denote the support of $\gamma$. For $\ell \geqslant 1$, we extend the notion of support to $\ell$-tuples of paths $\underline{\gamma} = (\gamma_1, \ldots, \gamma_\ell)$ as $\Supp(\underline{\gamma}) = \bigcup_{i=1}^\ell \Supp(\gamma_i)$. For $\ell \geqslant 1$, write $\Pi_\ell(G)$ the set of $\ell$-tuples of paths in $G$. 

From now on, assume that $G$ is acyclic, meaning there is no nontrivial path $\gamma$ in $G$ such that $s(\gamma) = t(\gamma)$. An \new{antichain} of $G$ is any subset of vertices $\{w_1,\ldots,w_r\} \subset G_0
$  such that there is no path $\gamma$ in $G$ with $s(\gamma) = w_i $ and $t(\gamma) = w_j$ for all $1 \leqslant i, j \leqslant r$ with $i \neq j$. 

A \new{filling} of $G$ is a map $f : G_0 \longrightarrow \mathbb{N}$. We assign to any $\ell$-tuple of paths $\underline{\gamma}$ of $G$ a \new{$f$-weight} defined by \[\wt_f(\underline{\gamma}) = \sum_{v \in \Supp(\underline{\gamma})} f(v).\]
Set $M_0^G(f) = 0$, and for all integers $\ell \geqslant 1$, $M_\ell^G(f) = \max_{\underline{\gamma} \in \Pi_\ell(G)} \wt_f(\underline{\gamma})$. We define the \new{Greene--Kleitman invariant} of $f$ in $G$ as \[\GK_G(f) = \left(M_\ell^G(f) - M_{\ell-1}^G(f) \right)_{\ell \geqslant 1}.\]

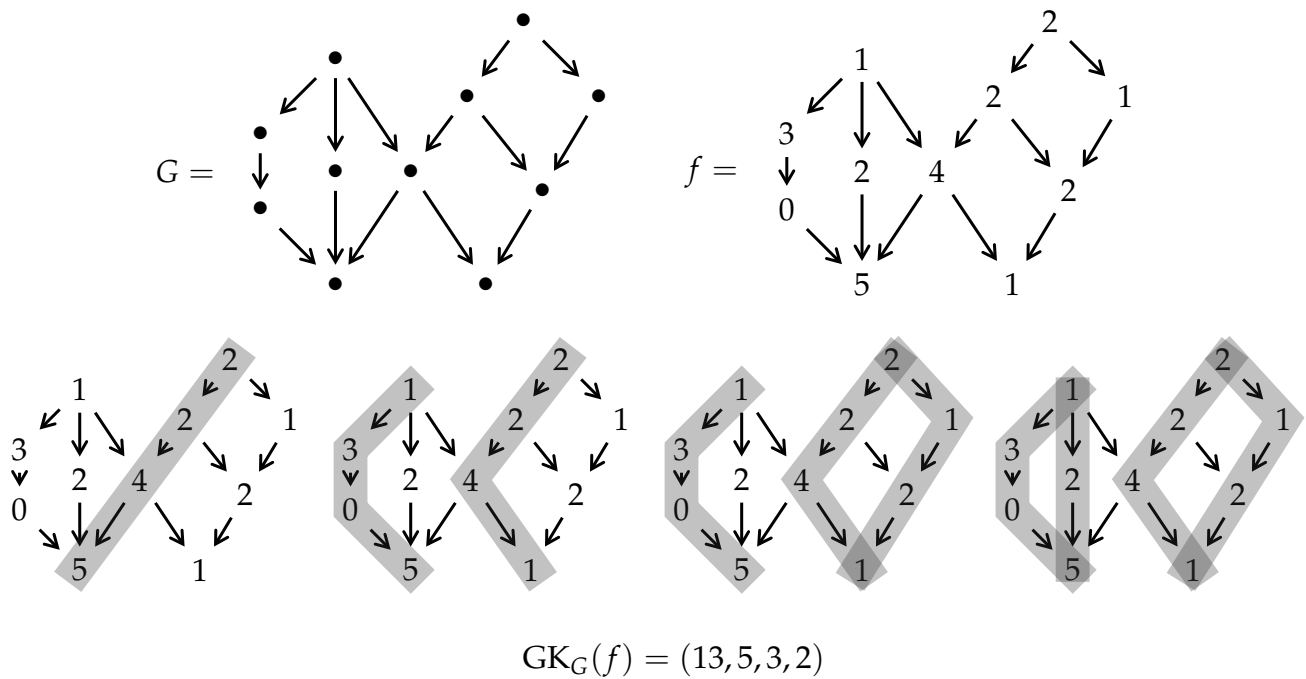
\begin{figure}
	\centering
	\begin{tikzpicture}[line width=0.4mm,>= angle 60]
		\node(0) at (-3,-1.5){$G =$};
		\node(1) at (-1,0){$\bullet$};
		\node(2) at (1.5,0.5){$\bullet$};
		\node(3) at (-2,-1){$\bullet$};
		\node(4) at (-1,-1.5){$\bullet$};
		\node(5) at (0.75,-0.5){$\bullet$};
		\node(6) at (2.5,-0.5){$\bullet$};
		\node(7) at (-2,-2){$\bullet$};
		\node(8) at (0,-1.5){$\bullet$};
		\node(9) at (1.75,-1.75){$\bullet$};
		\node(10) at (-1,-3){$\bullet$};
		\node(11) at (1,-3){$\bullet$};
		\draw[->] (1) -- (3);
		\draw[->] (1) -- (4);
		\draw[->] (1) -- (8);
		\draw[->] (2) -- (5);
		\draw[->] (2) -- (6);
		\draw[->] (3) -- (7);
		\draw[->] (4) -- (10);
		\draw[->] (5) -- (9);
		\draw[->] (6) -- (9);
		\draw[->] (7) -- (10);
		\draw[->] (8) -- (10);
		\draw[->] (8) -- (11);
		\draw[->] (9) -- (11);
		\draw[->] (5) -- (8);
		\begin{scope}[xshift=7cm]
			\node(1) at (-3,-1.5){$f =$};
			\node(1) at (-1,0){$1$};
			\node(2) at (1.5,0.5){$2$};
			\node(3) at (-2,-1){$3$};
			\node(4) at (-1,-1.5){$2$};
			\node(5) at (0.75,-0.5){$2$};
			\node(6) at (2.5,-0.5){$1$};
			\node(7) at (-2,-2){$0$};
			\node(8) at (0,-1.5){$4$};
			\node(9) at (1.75,-1.75){$2$};
			\node(10) at (-1,-3){$5$};
			\node(11) at (1,-3){$1$};
			\draw[->] (1) -- (3);
			\draw[->] (1) -- (4);
			\draw[->] (1) -- (8);
			\draw[->] (2) -- (5);
			\draw[->] (2) -- (6);
			\draw[->] (3) -- (7);
			\draw[->] (4) -- (10);
			\draw[->] (5) -- (9);
			\draw[->] (6) -- (9);
			\draw[->] (7) -- (10);
			\draw[->] (8) -- (10);
			\draw[->] (8) -- (11);
			\draw[->] (9) -- (11);
			\draw[->] (5) -- (8);
		\end{scope}
		
		\begin{scope}[scale = 0.8, xshift=-4.5cm, yshift=-5.5cm]
			\node(1) at (-1,0){$1$};
			\node(2) at (1.5,0.5){$2$};
			\node(3) at (-2,-1){$3$};
			\node(4) at (-1,-1.5){$2$};
			\node(5) at (0.75,-0.5){$2$};
			\node(6) at (2.5,-0.5){$1$};
			\node(7) at (-2,-2){$0$};
			\node(8) at (0,-1.5){$4$};
			\node(9) at (1.75,-1.75){$2$};
			\node(10) at (-1,-3){$5$};
			\node(11) at (1,-3){$1$};
			\draw[->] (1) -- (3);
			\draw[->] (1) -- (4);
			\draw[->] (1) -- (8);
			\draw[->] (2) -- (5);
			\draw[->] (2) -- (6);
			\draw[->] (3) -- (7);
			\draw[->] (4) -- (10);
			\draw[->] (5) -- (9);
			\draw[->] (6) -- (9);
			\draw[->] (7) -- (10);
			\draw[->] (8) -- (10);
			\draw[->] (8) -- (11);
			\draw[->] (9) -- (11);
			\draw[->] (5) -- (8);
			\draw[-,line width=4.5mm, black!80, draw opacity = 0.3] (1.7,0.7) -- (-1.2,-3.2);
		\end{scope}
		
		\begin{scope}[scale = 0.8, xshift=1cm, yshift=-5.5cm]
			\node(1) at (-1,0){$1$};
			\node(2) at (1.5,0.5){$2$};
			\node(3) at (-2,-1){$3$};
			\node(4) at (-1,-1.5){$2$};
			\node(5) at (0.75,-0.5){$2$};
			\node(6) at (2.5,-0.5){$1$};
			\node(7) at (-2,-2){$0$};
			\node(8) at (0,-1.5){$4$};
			\node(9) at (1.75,-1.75){$2$};
			\node(10) at (-1,-3){$5$};
			\node(11) at (1,-3){$1$};
			\draw[->] (1) -- (3);
			\draw[->] (1) -- (4);
			\draw[->] (1) -- (8);
			\draw[->] (2) -- (5);
			\draw[->] (2) -- (6);
			\draw[->] (3) -- (7);
			\draw[->] (4) -- (10);
			\draw[->] (5) -- (9);
			\draw[->] (6) -- (9);
			\draw[->] (7) -- (10);
			\draw[->] (8) -- (10);
			\draw[->] (8) -- (11);
			\draw[->] (9) -- (11);
			\draw[->] (5) -- (8);
			\draw[-,line width=4.5mm, black!80, draw opacity = 0.3] (1.7,0.7) -- (0,-1.5)-- (1.2,-3.2);
			\draw[-,line width=4.5mm, black!80, draw opacity = 0.3](-0.8,0.2) -- (-2,-1) -- (-2,-2) -- (-0.8,-3.2);
		\end{scope}
		
		\begin{scope}[scale = 0.8, xshift=6.5cm, yshift=-5.5cm]
			\node(1) at (-1,0){$1$};
			\node(2) at (1.5,0.5){$2$};
			\node(3) at (-2,-1){$3$};
			\node(4) at (-1,-1.5){$2$};
			\node(5) at (0.75,-0.5){$2$};
			\node(6) at (2.5,-0.5){$1$};
			\node(7) at (-2,-2){$0$};
			\node(8) at (0,-1.5){$4$};
			\node(9) at (1.75,-1.75){$2$};
			\node(10) at (-1,-3){$5$};
			\node(11) at (1,-3){$1$};
			\draw[->] (1) -- (3);
			\draw[->] (1) -- (4);
			\draw[->] (1) -- (8);
			\draw[->] (2) -- (5);
			\draw[->] (2) -- (6);
			\draw[->] (3) -- (7);
			\draw[->] (4) -- (10);
			\draw[->] (5) -- (9);
			\draw[->] (6) -- (9);
			\draw[->] (7) -- (10);
			\draw[->] (8) -- (10);
			\draw[->] (8) -- (11);
			\draw[->] (9) -- (11);
			\draw[->] (5) -- (8);
		\draw[-,line width=4.5mm, black!80, draw opacity = 0.3] (1.4,0.7) -- (2.5,-0.5) -- (0.8,-3.2);
		\draw[-,line width=4.5mm, black!80, draw opacity = 0.3] (1.7,0.7) -- (0,-1.5)-- (1.2,-3.2);
		\draw[-,line width=4.5mm, black!80, draw opacity = 0.3](-0.8,0.2) -- (-2,-1) -- (-2,-2) -- (-0.8,-3.2); 
		\end{scope}
		
		\begin{scope}[scale = 0.8, xshift=12cm, yshift=-5.5cm]
			\node(1) at (-1,0){$1$};
			\node(2) at (1.5,0.5){$2$};
			\node(3) at (-2,-1){$3$};
			\node(4) at (-1,-1.5){$2$};
			\node(5) at (0.75,-0.5){$2$};
			\node(6) at (2.5,-0.5){$1$};
			\node(7) at (-2,-2){$0$};
			\node(8) at (0,-1.5){$4$};
			\node(9) at (1.75,-1.75){$2$};
			\node(10) at (-1,-3){$5$};
			\node(11) at (1,-3){$1$};
			\draw[->] (1) -- (3);
			\draw[->] (1) -- (4);
			\draw[->] (1) -- (8);
			\draw[->] (2) -- (5);
			\draw[->] (2) -- (6);
			\draw[->] (3) -- (7);
			\draw[->] (4) -- (10);
			\draw[->] (5) -- (9);
			\draw[->] (6) -- (9);
			\draw[->] (7) -- (10);
			\draw[->] (8) -- (10);
			\draw[->] (8) -- (11);
			\draw[->] (9) -- (11);
			\draw[->] (5) -- (8);
			\draw[-,line width=4.5mm, black!80, draw opacity = 0.3] (1.4,0.7) -- (2.5,-0.5) -- (0.8,-3.2);
			\draw[-,line width=4.5mm, black!80, draw opacity = 0.3] (1.7,0.7) -- (0,-1.5)-- (1.2,-3.2);
			\draw[-,line width=4.5mm, black!80, draw opacity = 0.3](-0.8,0.2) -- (-2,-1) -- (-2,-2) -- (-0.8,-3.2); \draw[-,line width=4.5mm, black!80, draw opacity = 0.3](-1,0.2) -- (-1,-3.2);
		\end{scope}
		
		\node(1) at (3.5,-8){$\GK_G(f) = (13,5,3,2)$};
	\end{tikzpicture}
	\caption{\label{fig:GK} An example of the computation of $\GK_G$.}
\end{figure}

See \cref{fig:GK} for an explicit computation example.

\begin{prop}[Greene--Kleitman \cite{GK76}] \label{GKprop} 
	Let $G$ be a finite direct acyclic graph and $f$ be a filling of $G$.  The integer sequence $\GK_G(f)$ is an integer partition of length the maximal cardinality of an antichain in $G$.
\end{prop}

\subsection{Ferrers diagram RSK}

Throughout this section, we highlight Gansner's generalized version of the RSK correspondence, which gives, for any nonzero integer partition $\lambda$, a bijection from fillings of shape $\lambda$ to reverse plane partitions of shape $\lambda$. 

Fix a nonzero integer partition $\lambda$. Let $G_\lambda$ be the oriented acyclic graph such that:
\begin{enumerate}[label = $\bullet$]
	\item its vertices are the elements of $\Fer(\lambda)$;
	
	\item its arrows are given by:
	\begin{enumerate}[label = $\bullet$]
		\item $(i,j) \longrightarrow (i+1,j)$ whenever $(i,j), (i+1,j) \in \Fer(\lambda)$;
		
		\item $(i,j) \longrightarrow (i,j+1)$ whenever $(i,j),(i,j+1) \in \Fer(\lambda)$. 
	\end{enumerate}
\end{enumerate}
 For all $m \in \mathbb{Z}$, write $D_m(\lambda) = \{(i,j) \in \Fer(\lambda) \mid i-j + \lambda_1 = m \}$ for the $m$th diagonal of $\lambda$. Note that $D_m(\lambda) \neq \varnothing$ for $1 \leqslant m \leqslant h_\lambda(1,1)$, where $h_\lambda(1,1) = \#\{(i,j) \in \Fer(\lambda) \mid i=1 \text{ or } j=1\}$ denotes the \emph{hook length of the box $(1,1)$ in $\lambda$}. 

For each value $1 \leqslant m \leqslant h_\lambda(1,1)$, consider $(u_m, v_m)$ the maximal element of $D_m(\lambda)$. Write $G_\lambda(m)$ for the full subgraph of $G_\lambda$ given by the poset ideal generated by $(u_m,v_m)$. Note that $G_\lambda(m)$ admits only one source $(1,1)$, and only one sink $(u_m,v_m)$.

We define $g = \RSK_\lambda(f)$ to be the filling of shape $\lambda$ defined by \[\forall m \in \{1, \ldots, h_\lambda(1,1)\},\ \forall (i,j) \in D_m(\lambda),\quad g(i,j) = \GK_{G_\lambda(m)}(f)_{u_m-i+1}.\]

\begin{figure}[h!]
	\centering
		\scalebox{0.6}{
		\begin{tikzpicture}[scale=1]
			
			\tkzDefPoint(0,0){a}
			\tkzDefPoint(0,1){b}
			\tkzDefPoint(1,1){c}
			\tkzDefPoint(1,0){d}
			\tkzDrawPolygon[line width = 0.7mm, color = black](a,b,c,d);
			
			\tkzDefPoint(1,0){a}
			\tkzDefPoint(1,1){b}
			\tkzDefPoint(2,1){c}
			\tkzDefPoint(2,0){d}
			\tkzDrawPolygon[line width = 0.7mm, color = black](a,b,c,d);
			
			\tkzDefPoint(2,0){a}
			\tkzDefPoint(2,1){b}
			\tkzDefPoint(3,1){c}
			\tkzDefPoint(3,0){d}
			\tkzDrawPolygon[line width = 0.7mm, color = black](a,b,c,d);
			
			\tkzDefPoint(3,0){a}
			\tkzDefPoint(3,1){b}
			\tkzDefPoint(4,1){c}
			\tkzDefPoint(4,0){d}
			\tkzDrawPolygon[line width = 0.7mm, color = black](a,b,c,d);
			
			\tkzDefPoint(4,1){a}
			\tkzDefPoint(4,0){b}
			\tkzDefPoint(5,0){c}
			\tkzDefPoint(5,1){d}
			\tkzDrawPolygon[line width = 0.7mm, color = black](a,b,c,d);
			
			\tkzDefPoint(0,0){a}
			\tkzDefPoint(0,-1){b}
			\tkzDefPoint(1,-1){c}
			\tkzDefPoint(1,0){d}
			\tkzDrawPolygon[line width = 0.7mm, color = black](a,b,c,d);
			
			\tkzDefPoint(1,0){a}
			\tkzDefPoint(1,-1){b}
			\tkzDefPoint(2,-1){c}
			\tkzDefPoint(2,0){d}
			\tkzDrawPolygon[line width = 0.7mm, color = black](a,b,c,d);
			
			\tkzDefPoint(2,0){a}
			\tkzDefPoint(2,-1){b}
			\tkzDefPoint(3,-1){c}
			\tkzDefPoint(3,0){d}
			\tkzDrawPolygon[line width = 0.7mm, color = black](a,b,c,d);
			
			\tkzDefPoint(0,-2){a}
			\tkzDefPoint(0,-1){b}
			\tkzDefPoint(1,-1){c}
			\tkzDefPoint(1,-2){d}
			\tkzDrawPolygon[line width = 0.7mm, color = black](a,b,c,d);
			
			\tkzDefPoint(1,-2){a}
			\tkzDefPoint(1,-1){b}
			\tkzDefPoint(2,-1){c}
			\tkzDefPoint(2,-2){d}
			\tkzDrawPolygon[line width = 0.7mm, color = black](a,b,c,d);
			
			\tkzDefPoint(2,-2){a}
			\tkzDefPoint(2,-1){b}
			\tkzDefPoint(3,-1){c}
			\tkzDefPoint(3,-2){d}
			\tkzDrawPolygon[line width = 0.7mm, color = black](a,b,c,d);
			
			\tkzDefPoint(0,-3){a}
			\tkzDefPoint(0,-2){b}
			\tkzDefPoint(1,-2){c}
			\tkzDefPoint(1,-3){d}
			\tkzDrawPolygon[line width = 0.7mm, color = black](a,b,c,d);
			
			\tkzDefPoint(1,-3){a}
			\tkzDefPoint(1,-2){b}
			\tkzDefPoint(2,-2){c}
			\tkzDefPoint(2,-3){d}
			\tkzDrawPolygon[line width = 0.7mm, color = black](a,b,c,d);
			
			\tkzLabelPoint[below](0.5,0.9){{\Huge $1$}};
			\tkzLabelPoint[below](1.5,0.9){{\Huge $2$}};
			\tkzLabelPoint[below](2.5,0.9){{\Huge $1$}};
			\tkzLabelPoint[below](3.5,0.9){{\Huge $0$}};
			\tkzLabelPoint[below](4.5,0.9){{\Huge $3$}};
			\tkzLabelPoint[below](0.5,-0.1){{\Huge $2$}};
			\tkzLabelPoint[below](1.5,-0.1){{\Huge $1$}};
			\tkzLabelPoint[below](2.5,-0.1){{\Huge $1$}};
			\tkzLabelPoint[below](0.5,-1.1){{\Huge $2$}};
			\tkzLabelPoint[below](1.5,-1.1){{\Huge $1$}};
			\tkzLabelPoint[below](2.5,-1.1){{\Huge $3$}};
			\tkzLabelPoint[below](0.5,-2.1){{\Huge $3$}};
			\tkzLabelPoint[below](1.5,-2.1){{\Huge $2$}};
			
			\begin{scope}[xshift=9cm]
				
				\tkzDefPoint(0,0){a}
				\tkzDefPoint(0,1){b}
				\tkzDefPoint(1,1){c}
				\tkzDefPoint(1,0){d}
				\tkzDrawPolygon[line width = 0.7mm, color = black](a,b,c,d);
				
				\tkzDefPoint(1,0){a}
				\tkzDefPoint(1,1){b}
				\tkzDefPoint(2,1){c}
				\tkzDefPoint(2,0){d}
				\tkzDrawPolygon[line width = 0.7mm, color = black](a,b,c,d);
				
				\tkzDefPoint(2,0){a}
				\tkzDefPoint(2,1){b}
				\tkzDefPoint(3,1){c}
				\tkzDefPoint(3,0){d}
				\tkzDrawPolygon[line width = 0.7mm, color = black](a,b,c,d);
				
				\tkzDefPoint(3,0){a}
				\tkzDefPoint(3,1){b}
				\tkzDefPoint(4,1){c}
				\tkzDefPoint(4,0){d}
				\tkzDrawPolygon[line width = 0.7mm, color = black](a,b,c,d);
				
				\tkzDefPoint(4,1){a}
				\tkzDefPoint(4,0){b}
				\tkzDefPoint(5,0){c}
				\tkzDefPoint(5,1){d}
				\tkzDrawPolygon[line width = 0.7mm, color = black](a,b,c,d);
				
				\tkzDefPoint(0,0){a}
				\tkzDefPoint(0,-1){b}
				\tkzDefPoint(1,-1){c}
				\tkzDefPoint(1,0){d}
				\tkzDrawPolygon[line width = 0.7mm, color = black](a,b,c,d);
				
				\tkzDefPoint(1,0){a}
				\tkzDefPoint(1,-1){b}
				\tkzDefPoint(2,-1){c}
				\tkzDefPoint(2,0){d}
				\tkzDrawPolygon[line width = 0.7mm, color = black](a,b,c,d);
				
				\tkzDefPoint(2,0){a}
				\tkzDefPoint(2,-1){b}
				\tkzDefPoint(3,-1){c}
				\tkzDefPoint(3,0){d}
				\tkzDrawPolygon[line width = 0.7mm, color = black](a,b,c,d);
				
				\tkzDefPoint(0,-2){a}
				\tkzDefPoint(0,-1){b}
				\tkzDefPoint(1,-1){c}
				\tkzDefPoint(1,-2){d}
				\tkzDrawPolygon[line width = 0.7mm, color = black](a,b,c,d);
				
				\tkzDefPoint(1,-2){a}
				\tkzDefPoint(1,-1){b}
				\tkzDefPoint(2,-1){c}
				\tkzDefPoint(2,-2){d}
				\tkzDrawPolygon[line width = 0.7mm, color = black](a,b,c,d);
				
				\tkzDefPoint(2,-2){a}
				\tkzDefPoint(2,-1){b}
				\tkzDefPoint(3,-1){c}
				\tkzDefPoint(3,-2){d}
				\tkzDrawPolygon[line width = 0.7mm, color = black](a,b,c,d);
				
				\tkzDefPoint(0,-3){a}
				\tkzDefPoint(0,-2){b}
				\tkzDefPoint(1,-2){c}
				\tkzDefPoint(1,-3){d}
				\tkzDrawPolygon[line width = 0.7mm, color = black](a,b,c,d);
				
				\tkzDefPoint(1,-3){a}
				\tkzDefPoint(1,-2){b}
				\tkzDefPoint(2,-2){c}
				\tkzDefPoint(2,-3){d}
				\tkzDrawPolygon[line width = 0.7mm, color = black](a,b,c,d);
				
				\tkzLabelPoint[below](0.5,0.9){{\Huge $1$}};
				\tkzLabelPoint[below](1.5,0.9){{\Huge $3$}};
				\tkzLabelPoint[below](2.5,0.9){{\Huge $4$}};
				\tkzLabelPoint[below](3.5,0.9){{\Huge $4$}};
				\tkzLabelPoint[below](4.5,0.9){{\Huge $7$}};
				\tkzLabelPoint[below](0.5,-0.1){{\Huge $3$}};
				\tkzLabelPoint[below](1.5,-0.1){{\Huge $4$}};
				\tkzLabelPoint[below](2.5,-0.1){{\Huge $5$}};
				\tkzLabelPoint[below](0.5,-1.1){{\Huge $4$}};
				\tkzLabelPoint[below](1.5,-1.1){{\Huge $6$}};
				\tkzLabelPoint[below](2.5,-1.1){{\Huge $9$}};
				\tkzLabelPoint[below](0.5,-2.1){{\Huge $8$}};
				\tkzLabelPoint[below](1.5,-2.1){{\Huge $10$}};
				
			\end{scope}
			
			\draw [|->,line width=1.5mm,red] (5,-1) -- node[above]{{\Huge $\RSK_\lambda $}} (8,-1);
			
	\end{tikzpicture}}

$ $
	
	\scalebox{0.5}{
		\begin{tikzpicture}[scale=1]
			\begin{scope}[xshift=-1cm,yshift=6cm]
				\draw [line width=0.7mm, gray, dashed] (1.5,0.5) --  (3.5,-1.5);
				\node[gray] at (3.8,-1.8){{\Large $4$}};
				
				\draw [line width=0.7mm,gray, dashed] (2.5,0.5) --  (3.5,-0.5);
				\node[gray] at (3.8,-0.8){{\Large $3$}};
				
				\draw [line width=0.7mm,gray, dashed] (3.5,0.5) --  (4.5,-0.5);
				\node[gray] at (4.8,-0.8){{\Large $2$}};
				
				\draw [line width=0.7mm,red!70, dashed] (4.5,0.5) --  (5.5,-0.5);
				\node[red!70] at (5.8,-0.8){{\Large $1$}};
				
				\draw [line width=0.7mm, gray, dashed] (0.5,0.5) --  (3.5,-2.5);
				\node[gray] at (3.8,-2.8){{\Large $5$}};
				
				\draw [line width=0.7mm,gray, dashed] (0.5,-0.5) --  (2.5,-2.5);
				\node[gray] at (2.8,-2.8){{\Large $6$}};
				
				\draw [line width=0.7mm,gray, dashed] (0.5,-1.5) --  (2.5,-3.5);
				\node[gray] at (2.8,-3.8){{\Large $7$}};
				
				\draw [line width=0.7mm,gray, dashed] (0.5,-2.5) --  (1.5,-3.5);
				\node[gray] at (1.8,-3.8){{\Large $8$}};
				
				\tkzDefPoint(0,0){a}
				\tkzDefPoint(0,1){b}
				\tkzDefPoint(1,1){c}
				\tkzDefPoint(1,0){d}
				\tkzDrawPolygon[line width = 0.7mm, color = black](a,b,c,d);
				
				\tkzDefPoint(1,0){a}
				\tkzDefPoint(1,1){b}
				\tkzDefPoint(2,1){c}
				\tkzDefPoint(2,0){d}
				\tkzDrawPolygon[line width = 0.7mm, color = black](a,b,c,d);
				
				\tkzDefPoint(2,0){a}
				\tkzDefPoint(2,1){b}
				\tkzDefPoint(3,1){c}
				\tkzDefPoint(3,0){d}
				\tkzDrawPolygon[line width = 0.7mm, color = black](a,b,c,d);
				
				\tkzDefPoint(3,0){a}
				\tkzDefPoint(3,1){b}
				\tkzDefPoint(4,1){c}
				\tkzDefPoint(4,0){d}
				\tkzDrawPolygon[line width = 0.7mm, color = black](a,b,c,d);
				
				\tkzDefPoint(4,1){a}
				\tkzDefPoint(4,0){b}
				\tkzDefPoint(5,0){c}
				\tkzDefPoint(5,1){d}
				\tkzDrawPolygon[line width = 0.7mm, color = black, fill=red!70](a,b,c,d);
				
				\tkzDefPoint(0,0){a}
				\tkzDefPoint(0,-1){b}
				\tkzDefPoint(1,-1){c}
				\tkzDefPoint(1,0){d}
				\tkzDrawPolygon[line width = 0.7mm, color = black](a,b,c,d);
				
				\tkzDefPoint(1,0){a}
				\tkzDefPoint(1,-1){b}
				\tkzDefPoint(2,-1){c}
				\tkzDefPoint(2,0){d}
				\tkzDrawPolygon[line width = 0.7mm, color = black](a,b,c,d);
				
				\tkzDefPoint(2,0){a}
				\tkzDefPoint(2,-1){b}
				\tkzDefPoint(3,-1){c}
				\tkzDefPoint(3,0){d}
				\tkzDrawPolygon[line width = 0.7mm, color = black](a,b,c,d);
				
				\tkzDefPoint(0,-2){a}
				\tkzDefPoint(0,-1){b}
				\tkzDefPoint(1,-1){c}
				\tkzDefPoint(1,-2){d}
				\tkzDrawPolygon[line width = 0.7mm, color = black](a,b,c,d);
				
				\tkzDefPoint(1,-2){a}
				\tkzDefPoint(1,-1){b}
				\tkzDefPoint(2,-1){c}
				\tkzDefPoint(2,-2){d}
				\tkzDrawPolygon[line width = 0.7mm, color = black](a,b,c,d);
				
				\tkzDefPoint(2,-2){a}
				\tkzDefPoint(2,-1){b}
				\tkzDefPoint(3,-1){c}
				\tkzDefPoint(3,-2){d}
				\tkzDrawPolygon[line width = 0.7mm, color = black](a,b,c,d);
				
				\tkzDefPoint(0,-3){a}
				\tkzDefPoint(0,-2){b}
				\tkzDefPoint(1,-2){c}
				\tkzDefPoint(1,-3){d}
				\tkzDrawPolygon[line width = 0.7mm, color = black](a,b,c,d);
				
				\tkzDefPoint(1,-3){a}
				\tkzDefPoint(1,-2){b}
				\tkzDefPoint(2,-2){c}
				\tkzDefPoint(2,-3){d}
				\tkzDrawPolygon[line width = 0.7mm, color = black](a,b,c,d);

				\tkzLabelPoint[below](4.5,0.9){{\Huge $\mathbf{7}$}};
			\end{scope};
			
			\begin{scope}[xshift=-1cm,yshift = -1.5cm]
				
				\draw [line width=0.7mm, gray, dashed] (1.5,0.5) --  (3.5,-1.5);
				\node[gray] at (3.8,-1.8){{\Large $4$}};
				
				\draw [line width=0.7mm,gray, dashed] (2.5,0.5) --  (3.5,-0.5);
				\node[gray] at (3.8,-0.8){{\Large $3$}};
				
				\draw [line width=0.7mm,red!70, dashed] (3.5,0.5) --  (4.5,-0.5);
				\node[red!70] at (4.8,-0.8){{\Large $2$}};
				
				\draw [line width=0.7mm,gray, dashed] (4.5,0.5) --  (5.5,-0.5);
				\node[gray] at (5.8,-0.8){{\Large $1$}};
				
				\draw [line width=0.7mm, gray, dashed] (0.5,0.5) --  (3.5,-2.5);
				\node[gray] at (3.8,-2.8){{\Large $5$}};
				
				\draw [line width=0.7mm,gray, dashed] (0.5,-0.5) --  (2.5,-2.5);
				\node[gray] at (2.8,-2.8){{\Large $6$}};
				
				\draw [line width=0.7mm,gray, dashed] (0.5,-1.5) --  (2.5,-3.5);
				\node[gray] at (2.8,-3.8){{\Large $7$}};
				
				\draw [line width=0.7mm,gray, dashed] (0.5,-2.5) --  (1.5,-3.5);
				\node[gray] at (1.8,-3.8){{\Large $8$}};
				
				\tkzDefPoint(0,0){a}
				\tkzDefPoint(0,1){b}
				\tkzDefPoint(1,1){c}
				\tkzDefPoint(1,0){d}
				\tkzDrawPolygon[line width = 0.7mm, color = black](a,b,c,d);
				
				\tkzDefPoint(1,0){a}
				\tkzDefPoint(1,1){b}
				\tkzDefPoint(2,1){c}
				\tkzDefPoint(2,0){d}
				\tkzDrawPolygon[line width = 0.7mm, color = black](a,b,c,d);
				
				\tkzDefPoint(2,0){a}
				\tkzDefPoint(2,1){b}
				\tkzDefPoint(3,1){c}
				\tkzDefPoint(3,0){d}
				\tkzDrawPolygon[line width = 0.7mm, color = black](a,b,c,d);
				
				\tkzDefPoint(3,0){a}
				\tkzDefPoint(3,1){b}
				\tkzDefPoint(4,1){c}
				\tkzDefPoint(4,0){d}
				\tkzDrawPolygon[line width = 0.7mm, color = black, fill=red!70](a,b,c,d);
				
				\tkzDefPoint(4,1){a}
				\tkzDefPoint(4,0){b}
				\tkzDefPoint(5,0){c}
				\tkzDefPoint(5,1){d}
				\tkzDrawPolygon[line width = 0.7mm, color = black](a,b,c,d);
				
				\tkzDefPoint(0,0){a}
				\tkzDefPoint(0,-1){b}
				\tkzDefPoint(1,-1){c}
				\tkzDefPoint(1,0){d}
				\tkzDrawPolygon[line width = 0.7mm, color = black](a,b,c,d);
				
				\tkzDefPoint(1,0){a}
				\tkzDefPoint(1,-1){b}
				\tkzDefPoint(2,-1){c}
				\tkzDefPoint(2,0){d}
				\tkzDrawPolygon[line width = 0.7mm, color = black](a,b,c,d);
				
				\tkzDefPoint(2,0){a}
				\tkzDefPoint(2,-1){b}
				\tkzDefPoint(3,-1){c}
				\tkzDefPoint(3,0){d}
				\tkzDrawPolygon[line width = 0.7mm, color = black](a,b,c,d);
				
				\tkzDefPoint(0,-2){a}
				\tkzDefPoint(0,-1){b}
				\tkzDefPoint(1,-1){c}
				\tkzDefPoint(1,-2){d}
				\tkzDrawPolygon[line width = 0.7mm, color = black](a,b,c,d);
				
				\tkzDefPoint(1,-2){a}
				\tkzDefPoint(1,-1){b}
				\tkzDefPoint(2,-1){c}
				\tkzDefPoint(2,-2){d}
				\tkzDrawPolygon[line width = 0.7mm, color = black](a,b,c,d);
				
				\tkzDefPoint(2,-2){a}
				\tkzDefPoint(2,-1){b}
				\tkzDefPoint(3,-1){c}
				\tkzDefPoint(3,-2){d}
				\tkzDrawPolygon[line width = 0.7mm, color = black](a,b,c,d);
				
				\tkzDefPoint(0,-3){a}
				\tkzDefPoint(0,-2){b}
				\tkzDefPoint(1,-2){c}
				\tkzDefPoint(1,-3){d}
				\tkzDrawPolygon[line width = 0.7mm, color = black](a,b,c,d);
				
				\tkzDefPoint(1,-3){a}
				\tkzDefPoint(1,-2){b}
				\tkzDefPoint(2,-2){c}
				\tkzDefPoint(2,-3){d}
				\tkzDrawPolygon[line width = 0.7mm, color = black](a,b,c,d);
				
				\tkzLabelPoint[below](3.5,0.9){{\Huge $\mathbf{4}$}};
				\tkzLabelPoint[below](4.5,0.9){{\Huge $7$}};
				
			\end{scope};
			
			\begin{scope}[xshift=-1cm,yshift = -9cm]
				
				\draw [line width=0.7mm, gray, dashed] (1.5,0.5) --  (3.5,-1.5);
				\node[gray] at (3.8,-1.8){{\Large $4$}};
				
				\draw [line width=0.7mm,red!70, dashed] (2.5,0.5) --  (3.5,-0.5);
				\node[red!70] at (3.8,-0.8){{\Large $3$}};
				
				\draw [line width=0.7mm,gray, dashed] (3.5,0.5) --  (4.5,-0.5);
				\node[gray] at (4.8,-0.8){{\Large $2$}};
				
				\draw [line width=0.7mm,gray, dashed] (4.5,0.5) --  (5.5,-0.5);
				\node[gray] at (5.8,-0.8){{\Large $1$}};
				
				\draw [line width=0.7mm, gray, dashed] (0.5,0.5) --  (3.5,-2.5);
				\node[gray] at (3.8,-2.8){{\Large $5$}};
				
				\draw [line width=0.7mm,gray, dashed] (0.5,-0.5) --  (2.5,-2.5);
				\node[gray] at (2.8,-2.8){{\Large $6$}};
				
				\draw [line width=0.7mm,gray, dashed] (0.5,-1.5) --  (2.5,-3.5);
				\node[gray] at (2.8,-3.8){{\Large $7$}};
				
				\draw [line width=0.7mm,gray, dashed] (0.5,-2.5) --  (1.5,-3.5);
				\node[gray] at (1.8,-3.8){{\Large $8$}};
				
				\tkzDefPoint(0,0){a}
				\tkzDefPoint(0,1){b}
				\tkzDefPoint(1,1){c}
				\tkzDefPoint(1,0){d}
				\tkzDrawPolygon[line width = 0.7mm, color = black](a,b,c,d);
				
				\tkzDefPoint(1,0){a}
				\tkzDefPoint(1,1){b}
				\tkzDefPoint(2,1){c}
				\tkzDefPoint(2,0){d}
				\tkzDrawPolygon[line width = 0.7mm, color = black](a,b,c,d);
				
				\tkzDefPoint(2,0){a}
				\tkzDefPoint(2,1){b}
				\tkzDefPoint(3,1){c}
				\tkzDefPoint(3,0){d}
				\tkzDrawPolygon[line width = 0.7mm, color = black,fill=red!70](a,b,c,d);
				
				\tkzDefPoint(3,0){a}
				\tkzDefPoint(3,1){b}
				\tkzDefPoint(4,1){c}
				\tkzDefPoint(4,0){d}
				\tkzDrawPolygon[line width = 0.7mm, color = black](a,b,c,d);
				
				\tkzDefPoint(4,1){a}
				\tkzDefPoint(4,0){b}
				\tkzDefPoint(5,0){c}
				\tkzDefPoint(5,1){d}
				\tkzDrawPolygon[line width = 0.7mm, color = black](a,b,c,d);
				
				\tkzDefPoint(0,0){a}
				\tkzDefPoint(0,-1){b}
				\tkzDefPoint(1,-1){c}
				\tkzDefPoint(1,0){d}
				\tkzDrawPolygon[line width = 0.7mm, color = black](a,b,c,d);
				
				\tkzDefPoint(1,0){a}
				\tkzDefPoint(1,-1){b}
				\tkzDefPoint(2,-1){c}
				\tkzDefPoint(2,0){d}
				\tkzDrawPolygon[line width = 0.7mm, color = black](a,b,c,d);
				
				\tkzDefPoint(2,0){a}
				\tkzDefPoint(2,-1){b}
				\tkzDefPoint(3,-1){c}
				\tkzDefPoint(3,0){d}
				\tkzDrawPolygon[line width = 0.7mm, color = black](a,b,c,d);
				
				\tkzDefPoint(0,-2){a}
				\tkzDefPoint(0,-1){b}
				\tkzDefPoint(1,-1){c}
				\tkzDefPoint(1,-2){d}
				\tkzDrawPolygon[line width = 0.7mm, color = black](a,b,c,d);
				
				\tkzDefPoint(1,-2){a}
				\tkzDefPoint(1,-1){b}
				\tkzDefPoint(2,-1){c}
				\tkzDefPoint(2,-2){d}
				\tkzDrawPolygon[line width = 0.7mm, color = black](a,b,c,d);
				
				\tkzDefPoint(2,-2){a}
				\tkzDefPoint(2,-1){b}
				\tkzDefPoint(3,-1){c}
				\tkzDefPoint(3,-2){d}
				\tkzDrawPolygon[line width = 0.7mm, color = black](a,b,c,d);
				
				\tkzDefPoint(0,-3){a}
				\tkzDefPoint(0,-2){b}
				\tkzDefPoint(1,-2){c}
				\tkzDefPoint(1,-3){d}
				\tkzDrawPolygon[line width = 0.7mm, color = black](a,b,c,d);
				
				\tkzDefPoint(1,-3){a}
				\tkzDefPoint(1,-2){b}
				\tkzDefPoint(2,-2){c}
				\tkzDefPoint(2,-3){d}
				\tkzDrawPolygon[line width = 0.7mm, color = black](a,b,c,d);
				
				\tkzLabelPoint[below](2.5,0.9){{\Huge $\mathbf{4}$}};
				\tkzLabelPoint[below](3.5,0.9){{\Huge $4$}};
				\tkzLabelPoint[below](4.5,0.9){{\Huge $7$}};
			\end{scope};

			\begin{scope}[xshift=-1cm,yshift=-16.5cm]
				\draw [line width=0.7mm, red!70, dashed] (1.5,0.5) --  (3.5,-1.5);
				\node[red!70] at (3.8,-1.8){{\Large $4$}};
				
				\draw [line width=0.7mm,gray, dashed] (2.5,0.5) --  (3.5,-0.5);
				\node[gray] at (3.8,-0.8){{\Large $3$}};
				
				\draw [line width=0.7mm,gray, dashed] (3.5,0.5) --  (4.5,-0.5);
				\node[gray] at (4.8,-0.8){{\Large $2$}};
				
				\draw [line width=0.7mm,gray, dashed] (4.5,0.5) --  (5.5,-0.5);
				\node[gray] at (5.8,-0.8){{\Large $1$}};
				
				\draw [line width=0.7mm, gray, dashed] (0.5,0.5) --  (3.5,-2.5);
				\node[gray] at (3.8,-2.8){{\Large $5$}};
				
				\draw [line width=0.7mm,gray, dashed] (0.5,-0.5) --  (2.5,-2.5);
				\node[gray] at (2.8,-2.8){{\Large $6$}};
				
				\draw [line width=0.7mm,gray, dashed] (0.5,-1.5) --  (2.5,-3.5);
				\node[gray] at (2.8,-3.8){{\Large $7$}};
				
				\draw [line width=0.7mm,gray, dashed] (0.5,-2.5) --  (1.5,-3.5);
				\node[gray] at (1.8,-3.8){{\Large $8$}};
				
				\tkzDefPoint(0,0){a}
				\tkzDefPoint(0,1){b}
				\tkzDefPoint(1,1){c}
				\tkzDefPoint(1,0){d}
				\tkzDrawPolygon[line width = 0.7mm, color = black](a,b,c,d);
				
				\tkzDefPoint(1,0){a}
				\tkzDefPoint(1,1){b}
				\tkzDefPoint(2,1){c}
				\tkzDefPoint(2,0){d}
				\tkzDrawPolygon[line width = 0.7mm, color = black,fill=red!70](a,b,c,d);
				
				\tkzDefPoint(2,0){a}
				\tkzDefPoint(2,1){b}
				\tkzDefPoint(3,1){c}
				\tkzDefPoint(3,0){d}
				\tkzDrawPolygon[line width = 0.7mm, color = black](a,b,c,d);
				
				\tkzDefPoint(3,0){a}
				\tkzDefPoint(3,1){b}
				\tkzDefPoint(4,1){c}
				\tkzDefPoint(4,0){d}
				\tkzDrawPolygon[line width = 0.7mm, color = black](a,b,c,d);
				
				\tkzDefPoint(4,1){a}
				\tkzDefPoint(4,0){b}
				\tkzDefPoint(5,0){c}
				\tkzDefPoint(5,1){d}
				\tkzDrawPolygon[line width = 0.7mm, color = black](a,b,c,d);
				
				\tkzDefPoint(0,0){a}
				\tkzDefPoint(0,-1){b}
				\tkzDefPoint(1,-1){c}
				\tkzDefPoint(1,0){d}
				\tkzDrawPolygon[line width = 0.7mm, color = black](a,b,c,d);
				
				\tkzDefPoint(1,0){a}
				\tkzDefPoint(1,-1){b}
				\tkzDefPoint(2,-1){c}
				\tkzDefPoint(2,0){d}
				\tkzDrawPolygon[line width = 0.7mm, color = black](a,b,c,d);
				
				\tkzDefPoint(2,0){a}
				\tkzDefPoint(2,-1){b}
				\tkzDefPoint(3,-1){c}
				\tkzDefPoint(3,0){d}
				\tkzDrawPolygon[line width = 0.7mm, color = black,fill=red!70](a,b,c,d);
				
				\tkzDefPoint(0,-2){a}
				\tkzDefPoint(0,-1){b}
				\tkzDefPoint(1,-1){c}
				\tkzDefPoint(1,-2){d}
				\tkzDrawPolygon[line width = 0.7mm, color = black](a,b,c,d);
				
				\tkzDefPoint(1,-2){a}
				\tkzDefPoint(1,-1){b}
				\tkzDefPoint(2,-1){c}
				\tkzDefPoint(2,-2){d}
				\tkzDrawPolygon[line width = 0.7mm, color = black](a,b,c,d);
				
				\tkzDefPoint(2,-2){a}
				\tkzDefPoint(2,-1){b}
				\tkzDefPoint(3,-1){c}
				\tkzDefPoint(3,-2){d}
				\tkzDrawPolygon[line width = 0.7mm, color = black](a,b,c,d);
				
				\tkzDefPoint(0,-3){a}
				\tkzDefPoint(0,-2){b}
				\tkzDefPoint(1,-2){c}
				\tkzDefPoint(1,-3){d}
				\tkzDrawPolygon[line width = 0.7mm, color = black](a,b,c,d);
				
				\tkzDefPoint(1,-3){a}
				\tkzDefPoint(1,-2){b}
				\tkzDefPoint(2,-2){c}
				\tkzDefPoint(2,-3){d}
				\tkzDrawPolygon[line width = 0.7mm, color = black](a,b,c,d);
				
				\tkzLabelPoint[below](1.5,0.9){{\Huge $\mathbf{3}$}};
				\tkzLabelPoint[below](2.5,0.9){{\Huge $4$}};
				\tkzLabelPoint[below](3.5,0.9){{\Huge $4$}};
				\tkzLabelPoint[below](4.5,0.9){{\Huge $7$}};
				\tkzLabelPoint[below](2.5,-0.1){{\Huge $\mathbf{5}$}};
			\end{scope};

			\begin{scope}[yshift = 6cm,xshift=17cm]
				
				\draw [line width=0.7mm, gray, dashed] (1.5,0.5) --  (3.5,-1.5);
				\node[gray] at (3.8,-1.8){{\Large $4$}};
				
				\draw [line width=0.7mm,gray, dashed] (2.5,0.5) --  (3.5,-0.5);
				\node[gray] at (3.8,-0.8){{\Large $3$}};
				
				\draw [line width=0.7mm,gray, dashed] (3.5,0.5) --  (4.5,-0.5);
				\node[gray] at (4.8,-0.8){{\Large $2$}};
				
				\draw [line width=0.7mm,gray, dashed] (4.5,0.5) --  (5.5,-0.5);
				\node[gray] at (5.8,-0.8){{\Large $1$}};
				
				\draw [line width=0.7mm, red!70, dashed] (0.5,0.5) --  (3.5,-2.5);
				\node[red!70] at (3.8,-2.8){{\Large $5$}};
				
				\draw [line width=0.7mm,gray, dashed] (0.5,-0.5) --  (2.5,-2.5);
				\node[gray] at (2.8,-2.8){{\Large $6$}};
				
				\draw [line width=0.7mm,gray, dashed] (0.5,-1.5) --  (2.5,-3.5);
				\node[gray] at (2.8,-3.8){{\Large $7$}};
				
				\draw [line width=0.7mm,gray, dashed] (0.5,-2.5) --  (1.5,-3.5);
				\node[gray] at (1.8,-3.8){{\Large $8$}};
				
				\tkzDefPoint(0,0){a}
				\tkzDefPoint(0,1){b}
				\tkzDefPoint(1,1){c}
				\tkzDefPoint(1,0){d}
				\tkzDrawPolygon[line width = 0.7mm, color = black,fill=red!70](a,b,c,d);
				
				\tkzDefPoint(1,0){a}
				\tkzDefPoint(1,1){b}
				\tkzDefPoint(2,1){c}
				\tkzDefPoint(2,0){d}
				\tkzDrawPolygon[line width = 0.7mm, color = black](a,b,c,d);
				
				\tkzDefPoint(2,0){a}
				\tkzDefPoint(2,1){b}
				\tkzDefPoint(3,1){c}
				\tkzDefPoint(3,0){d}
				\tkzDrawPolygon[line width = 0.7mm, color = black](a,b,c,d);
				
				\tkzDefPoint(3,0){a}
				\tkzDefPoint(3,1){b}
				\tkzDefPoint(4,1){c}
				\tkzDefPoint(4,0){d}
				\tkzDrawPolygon[line width = 0.7mm, color = black](a,b,c,d);
				
				\tkzDefPoint(4,1){a}
				\tkzDefPoint(4,0){b}
				\tkzDefPoint(5,0){c}
				\tkzDefPoint(5,1){d}
				\tkzDrawPolygon[line width = 0.7mm, color = black](a,b,c,d);
				
				\tkzDefPoint(0,0){a}
				\tkzDefPoint(0,-1){b}
				\tkzDefPoint(1,-1){c}
				\tkzDefPoint(1,0){d}
				\tkzDrawPolygon[line width = 0.7mm, color = black](a,b,c,d);
				
				\tkzDefPoint(1,0){a}
				\tkzDefPoint(1,-1){b}
				\tkzDefPoint(2,-1){c}
				\tkzDefPoint(2,0){d}
				\tkzDrawPolygon[line width = 0.7mm, color = black,fill=red!70](a,b,c,d);
				
				\tkzDefPoint(2,0){a}
				\tkzDefPoint(2,-1){b}
				\tkzDefPoint(3,-1){c}
				\tkzDefPoint(3,0){d}
				\tkzDrawPolygon[line width = 0.7mm, color = black](a,b,c,d);
				
				\tkzDefPoint(0,-2){a}
				\tkzDefPoint(0,-1){b}
				\tkzDefPoint(1,-1){c}
				\tkzDefPoint(1,-2){d}
				\tkzDrawPolygon[line width = 0.7mm, color = black](a,b,c,d);
				
				\tkzDefPoint(1,-2){a}
				\tkzDefPoint(1,-1){b}
				\tkzDefPoint(2,-1){c}
				\tkzDefPoint(2,-2){d}
				\tkzDrawPolygon[line width = 0.7mm, color = black](a,b,c,d);
				
				\tkzDefPoint(2,-2){a}
				\tkzDefPoint(2,-1){b}
				\tkzDefPoint(3,-1){c}
				\tkzDefPoint(3,-2){d}
				\tkzDrawPolygon[line width = 0.7mm, color = black,fill=red!70](a,b,c,d);
				
				\tkzDefPoint(0,-3){a}
				\tkzDefPoint(0,-2){b}
				\tkzDefPoint(1,-2){c}
				\tkzDefPoint(1,-3){d}
				\tkzDrawPolygon[line width = 0.7mm, color = black](a,b,c,d);
				
				\tkzDefPoint(1,-3){a}
				\tkzDefPoint(1,-2){b}
				\tkzDefPoint(2,-2){c}
				\tkzDefPoint(2,-3){d}
				\tkzDrawPolygon[line width = 0.7mm, color = black](a,b,c,d);
				
				\tkzLabelPoint[below](0.5,0.9){{\Huge $\mathbf{1}$}};
				\tkzLabelPoint[below](1.5,0.9){{\Huge $3$}};
				\tkzLabelPoint[below](2.5,0.9){{\Huge $4$}};
				\tkzLabelPoint[below](3.5,0.9){{\Huge $4$}};
				\tkzLabelPoint[below](4.5,0.9){{\Huge $7$}};
				\tkzLabelPoint[below](1.5,-0.1){{\Huge $\mathbf{4}$}};
				\tkzLabelPoint[below](2.5,-0.1){{\Huge $5$}};
				\tkzLabelPoint[below](2.5,-1.1){{\Huge $\mathbf{9}$}};
				
			\end{scope};
			
			\begin{scope}[xshift=17cm, yshift = -1.5cm]
				\draw [line width=0.7mm, gray, dashed] (1.5,0.5) --  (3.5,-1.5);
				\node[gray] at (3.8,-1.8){{\Large $4$}};
				
				\draw [line width=0.7mm,gray, dashed] (2.5,0.5) --  (3.5,-0.5);
				\node[gray] at (3.8,-0.8){{\Large $3$}};
				
				\draw [line width=0.7mm,gray, dashed] (3.5,0.5) --  (4.5,-0.5);
				\node[gray] at (4.8,-0.8){{\Large $2$}};
				
				\draw [line width=0.7mm,gray, dashed] (4.5,0.5) --  (5.5,-0.5);
				\node[gray] at (5.8,-0.8){{\Large $1$}};
				
				\draw [line width=0.7mm, gray, dashed] (0.5,0.5) --  (3.5,-2.5);
				\node[gray] at (3.8,-2.8){{\Large $5$}};
				
				\draw [line width=0.7mm,red!70, dashed] (0.5,-0.5) --  (2.5,-2.5);
				\node[red!70] at (2.8,-2.8){{\Large $6$}};
				
				\draw [line width=0.7mm,gray, dashed] (0.5,-1.5) --  (2.5,-3.5);
				\node[gray] at (2.8,-3.8){{\Large $7$}};
				
				\draw [line width=0.7mm,gray, dashed] (0.5,-2.5) --  (1.5,-3.5);
				\node[gray] at (1.8,-3.8){{\Large $8$}};
				
				\tkzDefPoint(0,0){a}
				\tkzDefPoint(0,1){b}
				\tkzDefPoint(1,1){c}
				\tkzDefPoint(1,0){d}
				\tkzDrawPolygon[line width = 0.7mm, color = black](a,b,c,d);
				
				\tkzDefPoint(1,0){a}
				\tkzDefPoint(1,1){b}
				\tkzDefPoint(2,1){c}
				\tkzDefPoint(2,0){d}
				\tkzDrawPolygon[line width = 0.7mm, color = black](a,b,c,d);
				
				\tkzDefPoint(2,0){a}
				\tkzDefPoint(2,1){b}
				\tkzDefPoint(3,1){c}
				\tkzDefPoint(3,0){d}
				\tkzDrawPolygon[line width = 0.7mm, color = black](a,b,c,d);
				
				\tkzDefPoint(3,0){a}
				\tkzDefPoint(3,1){b}
				\tkzDefPoint(4,1){c}
				\tkzDefPoint(4,0){d}
				\tkzDrawPolygon[line width = 0.7mm, color = black](a,b,c,d);
				
				\tkzDefPoint(4,1){a}
				\tkzDefPoint(4,0){b}
				\tkzDefPoint(5,0){c}
				\tkzDefPoint(5,1){d}
				\tkzDrawPolygon[line width = 0.7mm, color = black](a,b,c,d);
				
				\tkzDefPoint(0,0){a}
				\tkzDefPoint(0,-1){b}
				\tkzDefPoint(1,-1){c}
				\tkzDefPoint(1,0){d}
				\tkzDrawPolygon[line width = 0.7mm, color = black,fill=red!70](a,b,c,d);
				
				\tkzDefPoint(1,0){a}
				\tkzDefPoint(1,-1){b}
				\tkzDefPoint(2,-1){c}
				\tkzDefPoint(2,0){d}
				\tkzDrawPolygon[line width = 0.7mm, color = black](a,b,c,d);
				
				\tkzDefPoint(2,0){a}
				\tkzDefPoint(2,-1){b}
				\tkzDefPoint(3,-1){c}
				\tkzDefPoint(3,0){d}
				\tkzDrawPolygon[line width = 0.7mm, color = black](a,b,c,d);
				
				\tkzDefPoint(0,-2){a}
				\tkzDefPoint(0,-1){b}
				\tkzDefPoint(1,-1){c}
				\tkzDefPoint(1,-2){d}
				\tkzDrawPolygon[line width = 0.7mm, color = black](a,b,c,d);
				
				\tkzDefPoint(1,-2){a}
				\tkzDefPoint(1,-1){b}
				\tkzDefPoint(2,-1){c}
				\tkzDefPoint(2,-2){d}
				\tkzDrawPolygon[line width = 0.7mm, color = black,fill=red!70](a,b,c,d);
				
				\tkzDefPoint(2,-2){a}
				\tkzDefPoint(2,-1){b}
				\tkzDefPoint(3,-1){c}
				\tkzDefPoint(3,-2){d}
				\tkzDrawPolygon[line width = 0.7mm, color = black](a,b,c,d);
				
				\tkzDefPoint(0,-3){a}
				\tkzDefPoint(0,-2){b}
				\tkzDefPoint(1,-2){c}
				\tkzDefPoint(1,-3){d}
				\tkzDrawPolygon[line width = 0.7mm, color = black](a,b,c,d);
				
				\tkzDefPoint(1,-3){a}
				\tkzDefPoint(1,-2){b}
				\tkzDefPoint(2,-2){c}
				\tkzDefPoint(2,-3){d}
				\tkzDrawPolygon[line width = 0.7mm, color = black](a,b,c,d);
				
				\tkzLabelPoint[below](0.5,0.9){{\Huge $1$}};
				\tkzLabelPoint[below](1.5,0.9){{\Huge $3$}};
				\tkzLabelPoint[below](2.5,0.9){{\Huge $4$}};
				\tkzLabelPoint[below](3.5,0.9){{\Huge $4$}};
				\tkzLabelPoint[below](4.5,0.9){{\Huge $7$}};
				\tkzLabelPoint[below](0.5,-0.1){{\Huge $\mathbf{3}$}};
				\tkzLabelPoint[below](1.5,-0.1){{\Huge $4$}};
				\tkzLabelPoint[below](2.5,-0.1){{\Huge $5$}};
				\tkzLabelPoint[below](1.5,-1.1){{\Huge $\mathbf{6}$}};
				\tkzLabelPoint[below](2.5,-1.1){{\Huge $9$}};
				
			\end{scope};
			
			\begin{scope}[xshift=17cm,yshift=-9cm]
				\draw [line width=0.7mm, gray, dashed] (1.5,0.5) --  (3.5,-1.5);
				\node[gray] at (3.8,-1.8){{\Large $4$}};
				
				\draw [line width=0.7mm,gray, dashed] (2.5,0.5) --  (3.5,-0.5);
				\node[gray] at (3.8,-0.8){{\Large $3$}};
				
				\draw [line width=0.7mm,gray, dashed] (3.5,0.5) --  (4.5,-0.5);
				\node[gray] at (4.8,-0.8){{\Large $2$}};
				
				\draw [line width=0.7mm,gray, dashed] (4.5,0.5) --  (5.5,-0.5);
				\node[gray] at (5.8,-0.8){{\Large $1$}};
				
				\draw [line width=0.7mm, gray, dashed] (0.5,0.5) --  (3.5,-2.5);
				\node[gray] at (3.8,-2.8){{\Large $5$}};
				
				\draw [line width=0.7mm,gray, dashed] (0.5,-0.5) --  (2.5,-2.5);
				\node[gray] at (2.8,-2.8){{\Large $6$}};
				
				\draw [line width=0.7mm,red!70, dashed] (0.5,-1.5) --  (2.5,-3.5);
				\node[red!70] at (2.8,-3.8){{\Large $7$}};
				
				\draw [line width=0.7mm,gray, dashed] (0.5,-2.5) --  (1.5,-3.5);
				\node[gray] at (1.8,-3.8){{\Large $8$}};
				
				\tkzDefPoint(0,0){a}
				\tkzDefPoint(0,1){b}
				\tkzDefPoint(1,1){c}
				\tkzDefPoint(1,0){d}
				\tkzDrawPolygon[line width = 0.7mm, color = black](a,b,c,d);
				
				\tkzDefPoint(1,0){a}
				\tkzDefPoint(1,1){b}
				\tkzDefPoint(2,1){c}
				\tkzDefPoint(2,0){d}
				\tkzDrawPolygon[line width = 0.7mm, color = black](a,b,c,d);
				
				\tkzDefPoint(2,0){a}
				\tkzDefPoint(2,1){b}
				\tkzDefPoint(3,1){c}
				\tkzDefPoint(3,0){d}
				\tkzDrawPolygon[line width = 0.7mm, color = black](a,b,c,d);
				
				\tkzDefPoint(3,0){a}
				\tkzDefPoint(3,1){b}
				\tkzDefPoint(4,1){c}
				\tkzDefPoint(4,0){d}
				\tkzDrawPolygon[line width = 0.7mm, color = black](a,b,c,d);
				
				\tkzDefPoint(4,1){a}
				\tkzDefPoint(4,0){b}
				\tkzDefPoint(5,0){c}
				\tkzDefPoint(5,1){d}
				\tkzDrawPolygon[line width = 0.7mm, color = black](a,b,c,d);
				
				\tkzDefPoint(0,0){a}
				\tkzDefPoint(0,-1){b}
				\tkzDefPoint(1,-1){c}
				\tkzDefPoint(1,0){d}
				\tkzDrawPolygon[line width = 0.7mm, color = black](a,b,c,d);
				
				\tkzDefPoint(1,0){a}
				\tkzDefPoint(1,-1){b}
				\tkzDefPoint(2,-1){c}
				\tkzDefPoint(2,0){d}
				\tkzDrawPolygon[line width = 0.7mm, color = black](a,b,c,d);
				
				\tkzDefPoint(2,0){a}
				\tkzDefPoint(2,-1){b}
				\tkzDefPoint(3,-1){c}
				\tkzDefPoint(3,0){d}
				\tkzDrawPolygon[line width = 0.7mm, color = black](a,b,c,d);
				
				\tkzDefPoint(0,-2){a}
				\tkzDefPoint(0,-1){b}
				\tkzDefPoint(1,-1){c}
				\tkzDefPoint(1,-2){d}
				\tkzDrawPolygon[line width = 0.7mm, color = black,fill=red!70](a,b,c,d);
				
				\tkzDefPoint(1,-2){a}
				\tkzDefPoint(1,-1){b}
				\tkzDefPoint(2,-1){c}
				\tkzDefPoint(2,-2){d}
				\tkzDrawPolygon[line width = 0.7mm, color = black](a,b,c,d);
				
				\tkzDefPoint(2,-2){a}
				\tkzDefPoint(2,-1){b}
				\tkzDefPoint(3,-1){c}
				\tkzDefPoint(3,-2){d}
				\tkzDrawPolygon[line width = 0.7mm, color = black](a,b,c,d);
				
				\tkzDefPoint(0,-3){a}
				\tkzDefPoint(0,-2){b}
				\tkzDefPoint(1,-2){c}
				\tkzDefPoint(1,-3){d}
				\tkzDrawPolygon[line width = 0.7mm, color = black](a,b,c,d);
				
				\tkzDefPoint(1,-3){a}
				\tkzDefPoint(1,-2){b}
				\tkzDefPoint(2,-2){c}
				\tkzDefPoint(2,-3){d}
				\tkzDrawPolygon[line width = 0.7mm, color = black,fill=red!70](a,b,c,d);
				
				\tkzLabelPoint[below](0.5,0.9){{\Huge $1$}};
				\tkzLabelPoint[below](1.5,0.9){{\Huge $3$}};
				\tkzLabelPoint[below](2.5,0.9){{\Huge $4$}};
				\tkzLabelPoint[below](3.5,0.9){{\Huge $4$}};
				\tkzLabelPoint[below](4.5,0.9){{\Huge $7$}};
				\tkzLabelPoint[below](0.5,-0.1){{\Huge $3$}};
				\tkzLabelPoint[below](1.5,-0.1){{\Huge $4$}};
				\tkzLabelPoint[below](2.5,-0.1){{\Huge $5$}};
				\tkzLabelPoint[below](0.5,-1.1){{\Huge $\mathbf{4}$}};
				\tkzLabelPoint[below](1.5,-1.1){{\Huge $6$}};
				\tkzLabelPoint[below](2.5,-1.1){{\Huge $9$}};
				\tkzLabelPoint[below](1.5,-2.1){{\Huge $\mathbf{10}$}};
			\end{scope}
			
			\begin{scope}[xshift=17cm,yshift=-16.5cm]
				\draw [line width=0.7mm, gray, dashed] (1.5,0.5) --  (3.5,-1.5);
				\node[gray] at (3.8,-1.8){{\Large $4$}};
				
				\draw [line width=0.7mm,gray, dashed] (2.5,0.5) --  (3.5,-0.5);
				\node[gray] at (3.8,-0.8){{\Large $3$}};
				
				\draw [line width=0.7mm,gray, dashed] (3.5,0.5) --  (4.5,-0.5);
				\node[gray] at (4.8,-0.8){{\Large $2$}};
				
				\draw [line width=0.7mm,gray, dashed] (4.5,0.5) --  (5.5,-0.5);
				\node[gray] at (5.8,-0.8){{\Large $1$}};
				
				\draw [line width=0.7mm, gray, dashed] (0.5,0.5) --  (3.5,-2.5);
				\node[gray] at (3.8,-2.8){{\Large $5$}};
				
				\draw [line width=0.7mm,gray, dashed] (0.5,-0.5) --  (2.5,-2.5);
				\node[gray] at (2.8,-2.8){{\Large $6$}};
				
				\draw [line width=0.7mm,gray, dashed] (0.5,-1.5) --  (2.5,-3.5);
				\node[gray] at (2.8,-3.8){{\Large $7$}};
				
				\draw [line width=0.7mm,red!70, dashed] (0.5,-2.5) --  (1.5,-3.5);
				\node[red!70] at (1.8,-3.8){{\Large $8$}};
				
				\tkzDefPoint(0,0){a}
				\tkzDefPoint(0,1){b}
				\tkzDefPoint(1,1){c}
				\tkzDefPoint(1,0){d}
				\tkzDrawPolygon[line width = 0.7mm, color = black](a,b,c,d);
				
				\tkzDefPoint(1,0){a}
				\tkzDefPoint(1,1){b}
				\tkzDefPoint(2,1){c}
				\tkzDefPoint(2,0){d}
				\tkzDrawPolygon[line width = 0.7mm, color = black](a,b,c,d);
				
				\tkzDefPoint(2,0){a}
				\tkzDefPoint(2,1){b}
				\tkzDefPoint(3,1){c}
				\tkzDefPoint(3,0){d}
				\tkzDrawPolygon[line width = 0.7mm, color = black](a,b,c,d);
				
				\tkzDefPoint(3,0){a}
				\tkzDefPoint(3,1){b}
				\tkzDefPoint(4,1){c}
				\tkzDefPoint(4,0){d}
				\tkzDrawPolygon[line width = 0.7mm, color = black](a,b,c,d);
				
				\tkzDefPoint(4,1){a}
				\tkzDefPoint(4,0){b}
				\tkzDefPoint(5,0){c}
				\tkzDefPoint(5,1){d}
				\tkzDrawPolygon[line width = 0.7mm, color = black](a,b,c,d);
				
				\tkzDefPoint(0,0){a}
				\tkzDefPoint(0,-1){b}
				\tkzDefPoint(1,-1){c}
				\tkzDefPoint(1,0){d}
				\tkzDrawPolygon[line width = 0.7mm, color = black](a,b,c,d);
				
				\tkzDefPoint(1,0){a}
				\tkzDefPoint(1,-1){b}
				\tkzDefPoint(2,-1){c}
				\tkzDefPoint(2,0){d}
				\tkzDrawPolygon[line width = 0.7mm, color = black](a,b,c,d);
				
				\tkzDefPoint(2,0){a}
				\tkzDefPoint(2,-1){b}
				\tkzDefPoint(3,-1){c}
				\tkzDefPoint(3,0){d}
				\tkzDrawPolygon[line width = 0.7mm, color = black](a,b,c,d);
				
				\tkzDefPoint(0,-2){a}
				\tkzDefPoint(0,-1){b}
				\tkzDefPoint(1,-1){c}
				\tkzDefPoint(1,-2){d}
				\tkzDrawPolygon[line width = 0.7mm, color = black](a,b,c,d);
				
				\tkzDefPoint(1,-2){a}
				\tkzDefPoint(1,-1){b}
				\tkzDefPoint(2,-1){c}
				\tkzDefPoint(2,-2){d}
				\tkzDrawPolygon[line width = 0.7mm, color = black](a,b,c,d);
				
				\tkzDefPoint(2,-2){a}
				\tkzDefPoint(2,-1){b}
				\tkzDefPoint(3,-1){c}
				\tkzDefPoint(3,-2){d}
				\tkzDrawPolygon[line width = 0.7mm, color = black](a,b,c,d);
				
				\tkzDefPoint(0,-3){a}
				\tkzDefPoint(0,-2){b}
				\tkzDefPoint(1,-2){c}
				\tkzDefPoint(1,-3){d}
				\tkzDrawPolygon[line width = 0.7mm, color = black,fill=red!70](a,b,c,d);
				
				\tkzDefPoint(1,-3){a}
				\tkzDefPoint(1,-2){b}
				\tkzDefPoint(2,-2){c}
				\tkzDefPoint(2,-3){d}
				\tkzDrawPolygon[line width = 0.7mm, color = black](a,b,c,d);
				
				\tkzLabelPoint[below](0.5,0.9){{\Huge $1$}};
				\tkzLabelPoint[below](1.5,0.9){{\Huge $3$}};
				\tkzLabelPoint[below](2.5,0.9){{\Huge $4$}};
				\tkzLabelPoint[below](3.5,0.9){{\Huge $4$}};
				\tkzLabelPoint[below](4.5,0.9){{\Huge $7$}};
				\tkzLabelPoint[below](0.5,-0.1){{\Huge $3$}};
				\tkzLabelPoint[below](1.5,-0.1){{\Huge $4$}};
				\tkzLabelPoint[below](2.5,-0.1){{\Huge $5$}};
				\tkzLabelPoint[below](0.5,-1.1){{\Huge $4$}};
				\tkzLabelPoint[below](1.5,-1.1){{\Huge $6$}};
				\tkzLabelPoint[below](2.5,-1.1){{\Huge $9$}};
				\tkzLabelPoint[below](0.5,-2.1){{\Huge $\mathbf{8}$}};
				\tkzLabelPoint[below](1.5,-2.1){{\Huge $10$}};
			\end{scope}
			
			\begin{scope}[->,line width=0.4mm,>= angle 60, scale =1.5,xshift=-6cm,yshift=5cm]
				
				\tkzDefPoint(-0.5,0.5){a}
				\tkzDefPoint(4.5,0.5){b}
				\tkzDefPoint(4.5,-0.5){c}
				\tkzDefPoint(-0.5,-0.5){d}
				\tkzDrawPolygon[line width = 2mm, color = red, fill = red!10](a,b,c,d);
				
				\node[circle,draw] (a) at (0,0){{\Huge $1$}};
				\node[circle,draw] (b) at (1,0){{\Huge $2$}};
				\node[circle,draw] (c) at (2,0){{\Huge $1$}};
				\node[circle,draw] (d) at (3,0){{\Huge $0$}};
				\node[circle,draw,fill=red!60] (e) at (4,0){{\Huge $3$}};
				
				\node[circle,draw](f) at (0,-1){{\Huge $2$}};
				\node[circle,draw] (g) at (1,-1){{\Huge $1$}};
				\node[circle,draw] (h) at (2,-1){{\Huge $1$}};
				
				\node[circle,draw] (i) at (0,-2){{\Huge $2$}};
				\node[circle,draw] (j) at (1,-2){{\Huge $1$}};
				\node[circle,draw] (k) at (2,-2){{\Huge $3$}};
				
				\node[circle,draw] (l) at (0,-3){{\Huge $3$}};
				\node[circle,draw] (m) at (1,-3){{\Huge $2$}};
				
				\draw (a) -- (b);
				\draw (a) -- (f);
				\draw (b) -- (c);
				\draw (b) -- (g);
				\draw (c) -- (d);
				\draw (c) -- (h);
				\draw (d) -- (e);
				\draw (f) -- (g);
				\draw (f) -- (i);
				\draw (g) -- (h);
				\draw (g) -- (j);
				\draw (h) -- (k);
				\draw (i) -- (l);
				\draw (i) -- (j);
				\draw (j) -- (m);
				\draw (j) -- (k);
				\draw (l) -- (m);
				\draw[-,line width=5mm, black!80, draw opacity = 0.3](0,0) edge (4,0);
			\end{scope};
			
			\begin{scope}[->,line width=0.4mm,>= angle 60, scale =1.5,xshift=-6cm,yshift=0cm]
				
				\tkzDefPoint(-0.5,0.5){a}
				\tkzDefPoint(3.5,0.5){b}
				\tkzDefPoint(3.5,-0.5){c}
				\tkzDefPoint(-0.5,-0.5){d}
				\tkzDrawPolygon[line width = 2mm, color = red, fill = red!10](a,b,c,d);
				
				\node[circle,draw] (a) at (0,0){{\Huge $1$}};
				\node[circle,draw] (b) at (1,0){{\Huge $2$}};
				\node[circle,draw] (c) at (2,0){{\Huge $1$}};
				\node[circle,draw,fill=red!60] (d) at (3,0){{\Huge $0$}};
				\node[circle,draw] (e) at (4,0){{\Huge $3$}};
				
				\node[circle,draw](f) at (0,-1){{\Huge $2$}};
				\node[circle,draw] (g) at (1,-1){{\Huge $1$}};
				\node[circle,draw] (h) at (2,-1){{\Huge $1$}};
				
				\node[circle,draw] (i) at (0,-2){{\Huge $2$}};
				\node[circle,draw] (j) at (1,-2){{\Huge $1$}};
				\node[circle,draw] (k) at (2,-2){{\Huge $3$}};
				
				\node[circle,draw] (l) at (0,-3){{\Huge $3$}};
				\node[circle,draw] (m) at (1,-3){{\Huge $2$}};
				
				\draw (a) -- (b);
				\draw (a) -- (f);
				\draw (b) -- (c);
				\draw (b) -- (g);
				\draw (c) -- (d);
				\draw (c) -- (h);
				\draw (d) -- (e);
				\draw (f) -- (g);
				\draw (f) -- (i);
				\draw (g) -- (h);
				\draw (g) -- (j);
				\draw (h) -- (k);
				\draw (i) -- (l);
				\draw (i) -- (j);
				\draw (j) -- (m);
				\draw (j) -- (k);
				\draw (l) -- (m);
				\draw[-,line width=5mm, black!80, draw opacity = 0.3](0,0) edge (3,0);
			\end{scope};
			
			\begin{scope}[->,line width=0.4mm,>= angle 60, scale =1.5,xshift=-6cm,yshift=-5cm]
				
				\tkzDefPoint(-0.5,0.5){a}
				\tkzDefPoint(2.5,0.5){b}
				\tkzDefPoint(2.5,-0.5){c}
				\tkzDefPoint(-0.5,-0.5){d}
				\tkzDrawPolygon[line width = 2mm, color = red, fill = red!10](a,b,c,d);
				
				\node[circle,draw] (a) at (0,0){{\Huge $1$}};
				\node[circle,draw] (b) at (1,0){{\Huge $2$}};
				\node[circle,draw,fill=red!60] (c) at (2,0){{\Huge $1$}};
				\node[circle,draw] (d) at (3,0){{\Huge $0$}};
				\node[circle,draw] (e) at (4,0){{\Huge $3$}};
				
				\node[circle,draw](f) at (0,-1){{\Huge $2$}};
				\node[circle,draw] (g) at (1,-1){{\Huge $1$}};
				\node[circle,draw] (h) at (2,-1){{\Huge $1$}};
				
				\node[circle,draw] (i) at (0,-2){{\Huge $2$}};
				\node[circle,draw] (j) at (1,-2){{\Huge $1$}};
				\node[circle,draw] (k) at (2,-2){{\Huge $3$}};
				
				\node[circle,draw] (l) at (0,-3){{\Huge $3$}};
				\node[circle,draw] (m) at (1,-3){{\Huge $2$}};
				
				\draw (a) -- (b);
				\draw (a) -- (f);
				\draw (b) -- (c);
				\draw (b) -- (g);
				\draw (c) -- (d);
				\draw (c) -- (h);
				\draw (d) -- (e);
				\draw (f) -- (g);
				\draw (f) -- (i);
				\draw (g) -- (h);
				\draw (g) -- (j);
				\draw (h) -- (k);
				\draw (i) -- (l);
				\draw (i) -- (j);
				\draw (j) -- (m);
				\draw (j) -- (k);
				\draw (l) -- (m);
				\draw[-,line width=5mm, black!80, draw opacity = 0.3](0,0) edge (2,0);
			\end{scope};
			
			\begin{scope}[->,line width=0.4mm,>= angle 60, scale =1.5,xshift=-6cm,yshift=-10cm]
				
				\tkzDefPoint(-0.5,0.5){a}
				\tkzDefPoint(2.5,0.5){b}
				\tkzDefPoint(2.5,-1.5){c}
				\tkzDefPoint(-0.5,-1.5){d}
				\tkzDrawPolygon[line width = 2mm, color = red, fill = red!10](a,b,c,d);
				
				\node[circle,draw] (a) at (0,0){{\Huge $1$}};
				\node[circle,draw] (b) at (1,0){{\Huge $2$}};
				\node[circle,draw] (c) at (2,0){{\Huge $1$}};
				\node[circle,draw] (d) at (3,0){{\Huge $0$}};
				\node[circle,draw] (e) at (4,0){{\Huge $3$}};
				
				\node[circle,draw](f) at (0,-1){{\Huge $2$}};
				\node[circle,draw] (g) at (1,-1){{\Huge $1$}};
				\node[circle,draw,fill=red!60] (h) at (2,-1){{\Huge $1$}};
				
				\node[circle,draw] (i) at (0,-2){{\Huge $2$}};
				\node[circle,draw] (j) at (1,-2){{\Huge $1$}};
				\node[circle,draw] (k) at (2,-2){{\Huge $3$}};
				
				\node[circle,draw] (l) at (0,-3){{\Huge $3$}};
				\node[circle,draw] (m) at (1,-3){{\Huge $2$}};
				
				\draw (a) -- (b);
				\draw (a) -- (f);
				\draw (b) -- (c);
				\draw (b) -- (g);
				\draw (c) -- (d);
				\draw (c) -- (h);
				\draw (d) -- (e);
				\draw (f) -- (g);
				\draw (f) -- (i);
				\draw (g) -- (h);
				\draw (g) -- (j);
				\draw (h) -- (k);
				\draw (i) -- (l);
				\draw (i) -- (j);
				\draw (j) -- (m);
				\draw (j) -- (k);
				\draw (l) -- (m);
				\draw[-,line width=5mm, black!80, draw opacity = 0.3](0,0) -- (2,0) -- (2,-1);
			\end{scope};
			
			\begin{scope}[->,line width=0.4mm,>= angle 60, scale =1.5,xshift=6cm,yshift=5cm]
				
				\tkzDefPoint(-0.5,0.5){a}
				\tkzDefPoint(2.5,0.5){b}
				\tkzDefPoint(2.5,-2.5){c}
				\tkzDefPoint(-0.5,-2.5){d}
				\tkzDrawPolygon[line width = 2mm, color = red, fill = red!10](a,b,c,d);
				
				\node[circle,draw] (a) at (0,0){{\Huge $1$}};
				\node[circle,draw] (b) at (1,0){{\Huge $2$}};
				\node[circle,draw] (c) at (2,0){{\Huge $1$}};
				\node[circle,draw] (d) at (3,0){{\Huge $0$}};
				\node[circle,draw] (e) at (4,0){{\Huge $3$}};
				
				\node[circle,draw](f) at (0,-1){{\Huge $2$}};
				\node[circle,draw] (g) at (1,-1){{\Huge $1$}};
				\node[circle,draw] (h) at (2,-1){{\Huge $1$}};
				
				\node[circle,draw] (i) at (0,-2){{\Huge $2$}};
				\node[circle,draw] (j) at (1,-2){{\Huge $1$}};
				\node[circle,draw,fill=red!60] (k) at (2,-2){{\Huge $3$}};
				
				\node[circle,draw] (l) at (0,-3){{\Huge $3$}};
				\node[circle,draw] (m) at (1,-3){{\Huge $2$}};
				
				\draw (a) -- (b);
				\draw (a) -- (f);
				\draw (b) -- (c);
				\draw (b) -- (g);
				\draw (c) -- (d);
				\draw (c) -- (h);
				\draw (d) -- (e);
				\draw (f) -- (g);
				\draw (f) -- (i);
				\draw (g) -- (h);
				\draw (g) -- (j);
				\draw (h) -- (k);
				\draw (i) -- (l);
				\draw (i) -- (j);
				\draw (j) -- (m);
				\draw (j) -- (k);
				\draw (l) -- (m);
				\draw[-,line width=5mm, black, draw opacity = 0.5](0,0) -- (0,-2) -- (2,-2);
				\draw[-,line width=5mm, black!80, draw opacity = 0.3](0,0) -- (2,0) -- (2,-2);
			\end{scope};
			
			\begin{scope}[->,line width=0.4mm,>= angle 60, scale =1.5,xshift=6cm,yshift=0cm]
				
				\tkzDefPoint(-0.5,0.5){a}
				\tkzDefPoint(1.5,0.5){b}
				\tkzDefPoint(1.5,-2.5){c}
				\tkzDefPoint(-0.5,-2.5){d}
				\tkzDrawPolygon[line width = 2mm, color = red, fill = red!10](a,b,c,d);
				
				\node[circle,draw] (a) at (0,0){{\Huge $1$}};
				\node[circle,draw] (b) at (1,0){{\Huge $2$}};
				\node[circle,draw] (c) at (2,0){{\Huge $1$}};
				\node[circle,draw] (d) at (3,0){{\Huge $0$}};
				\node[circle,draw] (e) at (4,0){{\Huge $3$}};
				
				\node[circle,draw](f) at (0,-1){{\Huge $2$}};
				\node[circle,draw] (g) at (1,-1){{\Huge $1$}};
				\node[circle,draw] (h) at (2,-1){{\Huge $1$}};
				
				\node[circle,draw] (i) at (0,-2){{\Huge $2$}};
				\node[circle,draw,fill=red!60] (j) at (1,-2){{\Huge $1$}};
				\node[circle,draw] (k) at (2,-2){{\Huge $3$}};
				
				\node[circle,draw] (l) at (0,-3){{\Huge $3$}};
				\node[circle,draw] (m) at (1,-3){{\Huge $2$}};
				
				\draw (a) -- (b);
				\draw (a) -- (f);
				\draw (b) -- (c);
				\draw (b) -- (g);
				\draw (c) -- (d);
				\draw (c) -- (h);
				\draw (d) -- (e);
				\draw (f) -- (g);
				\draw (f) -- (i);
				\draw (g) -- (h);
				\draw (g) -- (j);
				\draw (h) -- (k);
				\draw (i) -- (l);
				\draw (i) -- (j);
				\draw (j) -- (m);
				\draw (j) -- (k);
				\draw (l) -- (m);
				\draw[-,line width=5mm, black!80, draw opacity = 0.3](0,0) -- (0,-2) -- (1,-2);
			\end{scope};
			\begin{scope}[->,line width=0.4mm,>= angle 60, scale =1.5,xshift=6cm,yshift=-5cm]
				
				\tkzDefPoint(-0.5,0.5){a}
				\tkzDefPoint(1.5,0.5){b}
				\tkzDefPoint(1.5,-3.5){c}
				\tkzDefPoint(-0.5,-3.5){d}
				\tkzDrawPolygon[line width = 2mm, color = red, fill = red!10](a,b,c,d);
				
				\node[circle,draw] (a) at (0,0){{\Huge $1$}};
				\node[circle,draw] (b) at (1,0){{\Huge $2$}};
				\node[circle,draw] (c) at (2,0){{\Huge $1$}};
				\node[circle,draw] (d) at (3,0){{\Huge $0$}};
				\node[circle,draw] (e) at (4,0){{\Huge $3$}};
				
				\node[circle,draw](f) at (0,-1){{\Huge $2$}};
				\node[circle,draw] (g) at (1,-1){{\Huge $1$}};
				\node[circle,draw] (h) at (2,-1){{\Huge $1$}};
				
				\node[circle,draw] (i) at (0,-2){{\Huge $2$}};
				\node[circle,draw] (j) at (1,-2){{\Huge $1$}};
				\node[circle,draw] (k) at (2,-2){{\Huge $3$}};
				
				\node[circle,draw] (l) at (0,-3){{\Huge $3$}};
				\node[circle,draw,fill=red!60] (m) at (1,-3){{\Huge $2$}};
				
				\draw (a) -- (b);
				\draw (a) -- (f);
				\draw (b) -- (c);
				\draw (b) -- (g);
				\draw (c) -- (d);
				\draw (c) -- (h);
				\draw (d) -- (e);
				\draw (f) -- (g);
				\draw (f) -- (i);
				\draw (g) -- (h);
				\draw (g) -- (j);
				\draw (h) -- (k);
				\draw (i) -- (l);
				\draw (i) -- (j);
				\draw (j) -- (m);
				\draw (j) -- (k);
				\draw (l) -- (m);
				\draw[-,line width=5mm, black!80, draw opacity = 0.3](0,0) -- (0,-3) -- (1,-3);
			\end{scope};
			\begin{scope}[->,line width=0.4mm,>= angle 60, scale =1.5,xshift=6cm,yshift=-10cm]
				
				\tkzDefPoint(-0.5,0.5){a}
				\tkzDefPoint(0.5,0.5){b}
				\tkzDefPoint(0.5,-3.5){c}
				\tkzDefPoint(-0.5,-3.5){d}
				\tkzDrawPolygon[line width = 2mm, color = red, fill = red!10](a,b,c,d);
				
				\node[circle,draw] (a) at (0,0){{\Huge $1$}};
				\node[circle,draw] (b) at (1,0){{\Huge $2$}};
				\node[circle,draw] (c) at (2,0){{\Huge $1$}};
				\node[circle,draw] (d) at (3,0){{\Huge $0$}};
				\node[circle,draw] (e) at (4,0){{\Huge $3$}};
				
				\node[circle,draw](f) at (0,-1){{\Huge $2$}};
				\node[circle,draw] (g) at (1,-1){{\Huge $1$}};
				\node[circle,draw] (h) at (2,-1){{\Huge $1$}};
				
				\node[circle,draw] (i) at (0,-2){{\Huge $2$}};
				\node[circle,draw] (j) at (1,-2){{\Huge $1$}};
				\node[circle,draw] (k) at (2,-2){{\Huge $3$}};
				
				\node[circle,draw,fill=red!60] (l) at (0,-3){{\Huge $3$}};
				\node[circle,draw] (m) at (1,-3){{\Huge $2$}};
				
				\draw (a) -- (b);
				\draw (a) -- (f);
				\draw (b) -- (c);
				\draw (b) -- (g);
				\draw (c) -- (d);
				\draw (c) -- (h);
				\draw (d) -- (e);
				\draw (f) -- (g);
				\draw (f) -- (i);
				\draw (g) -- (h);
				\draw (g) -- (j);
				\draw (h) -- (k);
				\draw (i) -- (l);
				\draw (i) -- (j);
				\draw (j) -- (m);
				\draw (j) -- (k);
				\draw (l) -- (m);
				\draw[-,line width=5mm, black!80, draw opacity = 0.3](0,0) -- (0,-3);
			\end{scope};
	\end{tikzpicture}}
\caption{\label{fig:genRSK} Explicit calculations of $\RSK_\lambda(f)$ for a given filling $f$ of shape $\lambda = (5,3,3,2)$. For $1 \leqslant m \leqslant 8$, each framed subgraph corresponds to the subgraph $G_\lambda(m)$, and each filled diagonal colored in red corresponds to $\GK_{G_\lambda(m)}(f)$.} 
\end{figure}
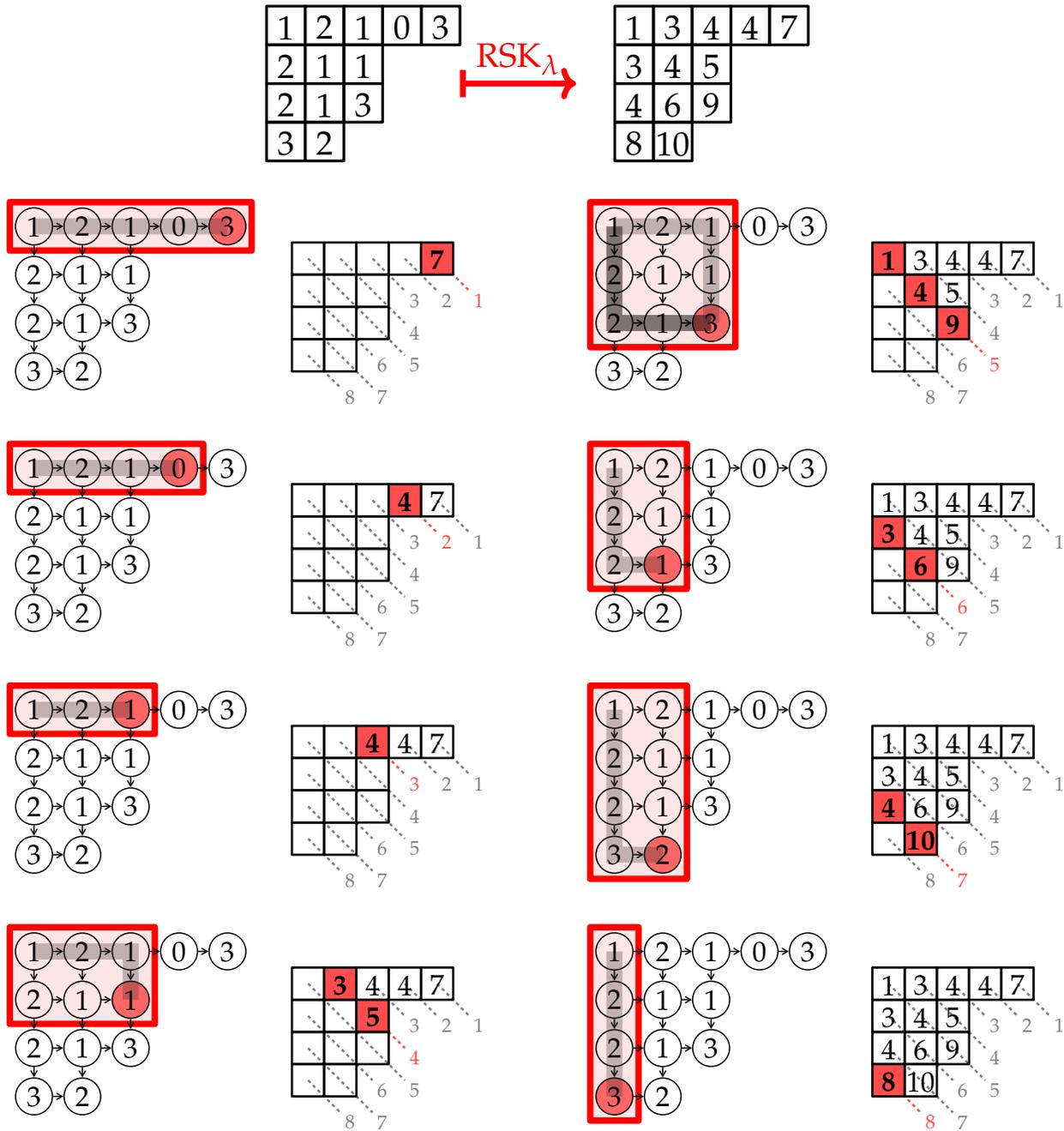

See \cref{fig:genRSK} for an explicit calculation of $\RSK_\lambda(f)$ for a given filling of $\lambda = (5,3,3,2)$.

\begin{thm}[Gansner \cite{Ga81Hi}]
	Let $\lambda$ be a nonzero integer partition. The map $\RSK_\lambda$ is a bijection from fillings of shape $\lambda$ to reverse plane partitions of shape $\lambda$.
\end{thm}

\begin{rmk}
	If $\lambda$ is a rectangle, we can recover the classical RSK. See \cite{GK76} and \cite[Section 6]{GPT19} for more details.
	
	Moreover, a parallel can be made with Britz and Fomin's version of the RSK algorithm \cite{BF99}, where we compute sequences of integer partitions for an $n \times n$ nonnegative integer matrix as growth diagrams. A generalized version of $\RSK$ was also exploited by Krattenthaler \cite{Kr06} on polyominos. From a given filling $f$ of shape $\lambda$, the integer partitions we can read on diagonals $D_m(\lambda)$ of $\RSK_\lambda(f)$ correspond precisely to the results obtained at the end of each line by using the Krattenthaler growth diagram algorithm version.
\end{rmk}

\section{Some tools}

In this section, we give the definition of some combinatorial objects that will be useful to present our generalized version of Gansner's RSK correspondence.

\subsection{Interval bipartitions}

An \new{interval bipartition} is a pair $(\B,\E) \in \mathcal{P}(\mathbb{N}^*)^2$ such that $\{\B,\E\}$ is a set partition of $\{i,\ldots,j\}$ for some $1 \leqslant i \leqslant j$. Call it \new{elementary} whenever $1 \in \B$ and $\max(\B \cup \E) \in \E$.

\begin{figure}[h!]
	\centering
		\scalebox{0.6}{\begin{tikzpicture}
				\tkzDefPoint(0,0){a}
				\tkzDefPoint(0,1){b}
				\tkzDefPoint(1,1){c}
				\tkzDefPoint(1,0){d}
				\tkzDrawPolygon[line width = 0.7mm, color = black](a,b,c,d);
				
				\tkzDefPoint(1,0){a}
				\tkzDefPoint(1,1){b}
				\tkzDefPoint(2,1){c}
				\tkzDefPoint(2,0){d}
				\tkzDrawPolygon[line width = 0.7mm, color = black](a,b,c,d);
				
				\tkzDefPoint(2,0){a}
				\tkzDefPoint(2,1){b}
				\tkzDefPoint(3,1){c}
				\tkzDefPoint(3,0){d}
				\tkzDrawPolygon[line width = 0.7mm, color = black](a,b,c,d);
				
				\tkzDefPoint(3,0){a}
				\tkzDefPoint(3,1){b}
				\tkzDefPoint(4,1){c}
				\tkzDefPoint(4,0){d}
				\tkzDrawPolygon[line width = 0.7mm, color = black](a,b,c,d);
				
				\tkzDefPoint(4,1){a}
				\tkzDefPoint(4,0){b}
				\tkzDefPoint(5,0){c}
				\tkzDefPoint(5,1){d}
				\tkzDrawPolygon[line width = 0.7mm, color = black](a,b,c,d);
				
				\tkzDefPoint(0,0){a}
				\tkzDefPoint(0,-1){b}
				\tkzDefPoint(1,-1){c}
				\tkzDefPoint(1,0){d}
				\tkzDrawPolygon[line width = 0.7mm, color = black](a,b,c,d);
				
				\tkzDefPoint(1,0){a}
				\tkzDefPoint(1,-1){b}
				\tkzDefPoint(2,-1){c}
				\tkzDefPoint(2,0){d}
				\tkzDrawPolygon[line width = 0.7mm, color = black](a,b,c,d);
				
				\tkzDefPoint(2,0){a}
				\tkzDefPoint(2,-1){b}
				\tkzDefPoint(3,-1){c}
				\tkzDefPoint(3,0){d}
				\tkzDrawPolygon[line width = 0.7mm, color = black](a,b,c,d);
				
				\tkzDefPoint(3,0){a}
				\tkzDefPoint(3,-1){b}
				\tkzDefPoint(4,-1){c}
				\tkzDefPoint(4,0){d}
				\tkzDrawPolygon[line width = 0.7mm, color = black](a,b,c,d);
				
				\tkzDefPoint(4,-1){a}
				\tkzDefPoint(4,0){b}
				\tkzDefPoint(5,0){c}
				\tkzDefPoint(5,-1){d}
				\tkzDrawPolygon[line width = 0.7mm, color = black](a,b,c,d);
				
				\tkzDefPoint(0,-2){a}
				\tkzDefPoint(0,-1){b}
				\tkzDefPoint(1,-1){c}
				\tkzDefPoint(1,-2){d}
				\tkzDrawPolygon[line width = 0.7mm, color = black](a,b,c,d);
				
				\tkzDefPoint(1,-2){a}
				\tkzDefPoint(1,-1){b}
				\tkzDefPoint(2,-1){c}
				\tkzDefPoint(2,-2){d}
				\tkzDrawPolygon[line width = 0.7mm, color = black](a,b,c,d);
				
				\tkzDefPoint(2,-2){a}
				\tkzDefPoint(2,-1){b}
				\tkzDefPoint(3,-1){c}
				\tkzDefPoint(3,-2){d}
				\tkzDrawPolygon[line width = 0.7mm, color = black](a,b,c,d);
				
				\tkzDefPoint(3,-2){a}
				\tkzDefPoint(3,-1){b}
				\tkzDefPoint(4,-1){c}
				\tkzDefPoint(4,-2){d}
				\tkzDrawPolygon[line width = 0.7mm, color = black](a,b,c,d);
				
				\tkzDefPoint(0,-3){a}
				\tkzDefPoint(0,-2){b}
				\tkzDefPoint(1,-2){c}
				\tkzDefPoint(1,-3){d}
				\tkzDrawPolygon[line width = 0.7mm, color = black](a,b,c,d);
				
				\node at (0.5,1.5){{\Large $9$}};
				\node at (1.5,1.5){{\Large $7$}};
				\node at (2.5,1.5){{\Large $6$}};
				\node at (3.5,1.5){{\Large $5$}};
				\node at (4.5,1.5){{\Large $3$}};
				\node at (-0.5,0.5){{\Large $1$}};
				\node at (-0.5,-0.5){{\Large $2$}};
				\node at (-0.5,-1.5){{\Large $4$}};
				\node at (-0.5,-2.5){{\Large $8$}};
		\end{tikzpicture}}
	\caption{\label{fig:BElambda} The (labelled) integer partition $\lambda(\B,\E)$ with $\B = \{1,2,4,8\}$ and $\E = \{3,5,6,7,9\}$.}
\end{figure}
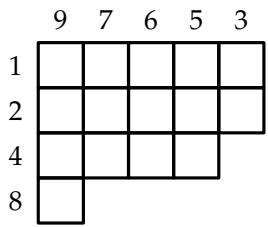

Fix $(\B,\E)$ as an interval bipartition. Write $\B = \{b_1 < b_2 < \ldots < b_p\}$.  We define the integer partition $\lambda(\B,\E)$ by $\lambda(\B,\E)_i = \#\{e \in \E\ \mid b_i < e\}$. If we also write $\E = \{e_1 < \ldots < e_q\}$, we can also describe $\lambda(\B,\E)$ by its Ferrers diagram: we have $(i,j) \in \Fer(\lambda(\B,\E))$ whenever $b_i < e_{q-j+1}$. It allows us to label the $i$th row of $\Fer(\lambda(\B,\E))$ by $b_i$ and the $j$th row by $e_{q-j+1}$. See \cref{fig:BElambda} for an example of such an object.

\begin{prop}
	For any integer partition $\lambda$, there exists an interval bipartition $(\B,\E)$ such that $\lambda(\B,\E) = \lambda$. Moreover, if $\lambda$ is a nonzero integer partition, there exists a unique elementary interval bipartition satisfying this property.
\end{prop}

\subsection{(Type $A$) Coxeter elements}

For any $n \geqslant 2$, let $\mathfrak{S}_n$ be the symmetric group on $n$ letters. For $1 \leqslant i < j \leqslant n$, write $(i,j)$ for the transposition exchanging $i$ and $j$. For $1 \leqslant i < n$, let $s_i$ be the adjacent transposition $(i,i+1)$. Let $S$ be the set of the adjacent transpositions. 

For any $w \in \mathfrak{S}_n$, an expression of $w$ is a way to write $w$ as a product of adjacent transpositions in $S$. The length $\ell(w)$ of $w$ is the minimal number of adjacent transpositions in $S$ needed to express $w$. Whenever, for some $1 \leqslant i < n$, $\ell(s_i w) < \ell(w)$, we say that $s_i$ is initial in $w$. Similarly, we call $s_i$ final in $w$ whenever $\ell(w s_i) < \ell(w)$.

A \new{Coxeter element (of $\mathfrak{S}_n$)} is an element $c \in \mathfrak{S}_n$ which can be written as a product of all the adjacent transpositions, in some order, where each of them appears exactly once. For example, $c = s_2 s_1 s_3 s_6 s_5 s_4 s_8 s_7 = (1,3,4,7,9,8,6,5,2)$ is a Coxeter element of $\mathfrak{S}_9$.

\begin{lem}
	An element $c \in \mathfrak{S}_n$ is a Coxeter element if and only if $c$ is a long cycle which can be written as follows \[c = (c_1,c_2, \ldots, c_m, c_{m+1}, \ldots, c_n)\]
	where $c_1 =1 < c_2 <\ldots <c_m = n > c_{m+1} > \ldots > c_n > c_1 = 1$.
\end{lem}

\subsection{Auslander--Reiten quivers}

\begin{figure}[h!]
	\centering
	\[
	\scalebox{0.9}{
		\begin{tikzpicture}

			\node (67) at (0,2) {(67)};
			\node (57) at (1,3) {(57)};
			\node (27) at (2,4) {(27)};
			\node (17) at (3,5) {(17)};
			\node (37) at (4,6) {(37)};
			\node (47) at (5,7) {(47)};
			
			\node (69) at (1,1) {(69)};
			\node (59) at (2,2) {(59)};
			\node (29) at (3,3) {(29)};
			\node (19) at (4,4) {(19)};
			\node (39) at (5,5) {(39)};
			
			\node (49) at (6,6) {(49)};
			
			\node (68) at (2,0) {(68)};
			\node (58) at (3,1) {(58)};
			\node (28) at (4,2) {(28)};
			\node (18) at (5,3) {(18)};
			\node (38) at (6,4) {(38)};
			\node (48) at (7,5) {(48)};
			\node (78) at (8,6) {(78)};
			\node (89) at (0,0) {(89)};
			\node (79) at (7,7) {(79)};
			
			\node (56) at (4,0) {(56)};
			\node (25) at (6,0) {(25)};
			\node (12) at (8,0) {(12)};
			\node (13) at (1,7) {(13)};
			\node (34) at (3,7) {(34)};
			
			\node (26) at (5,1) {(26)};
			\node (15) at (7,1) {(15)};
			\node (23) at (0,6) {(23)};
			\node (14) at (2,6) {(14)};
			
			\node (16) at (6,2) {(16)};
			\node (35) at (8,2) {(35)};
			\node (24) at (1,5) {(24)};
			
			\node (36) at (7,3) {(36)};
			\node (45) at (9,3) {(45)};
			
			\node (46) at (8,4) {(46)};
			
			\draw[->] (13)--(14);\draw[->] (14)--(34);\draw[->] (14)--(17);\draw[->] (15)--(35);\draw[->] (15)--(12);\draw[->] (16)--(36);\draw[->] (16)--(15);\draw[->] (17)--(37);\draw[->] (17)--(19);\draw[->] (18)--(38);\draw[->] (18)--(16);\draw[->] (19)--(39);\draw[->] (19)--(18);\draw[->] (23)--(13);\draw[->] (23)--(24);\draw[->] (24)--(14);\draw[->] (24)--(27);\draw[->] (25)--(15);\draw[->] (26)--(16);\draw[->] (26)--(25);\draw[->] (27)--(17);\draw[->] (27)--(29);\draw[->] (28)--(18);\draw[->] (28)--(26);\draw[->] (29)--(19);\draw[->] (29)--(28);\draw[->] (34)--(37);\draw[->] (35)--(45);\draw[->] (36)--(46);\draw[->] (36)--(35);\draw[->] (37)--(47);\draw[->] (37)--(39);\draw[->] (38)--(48);\draw[->] (38)--(36);\draw[->] (39)--(49);\draw[->] (39)--(38);\draw[->] (46)--(45);\draw[->] (47)--(49);\draw[->] (48)--(78);\draw[->] (48)--(46);\draw[->] (49)--(79);\draw[->] (49)--(48);\draw[->] (56)--(26);\draw[->] (57)--(27);\draw[->] (57)--(59);\draw[->] (58)--(28);\draw[->] (58)--(56);\draw[->] (59)--(29);\draw[->] (59)--(58);\draw[->] (67)--(57);\draw[->] (67)--(69);\draw[->] (68)--(58);\draw[->] (69)--(59);\draw[->] (69)--(68);\draw[->] (79)--(78);\draw[->] (89)--(69);
	\end{tikzpicture} }\]
\caption{\label{fig:ARc} The Auslander--Reiten quiver of $c = (1,3,4,7,9,8,6,5,2)=s_2s_1s_3s_6s_5s_4s_8s_7$.}
\end{figure}
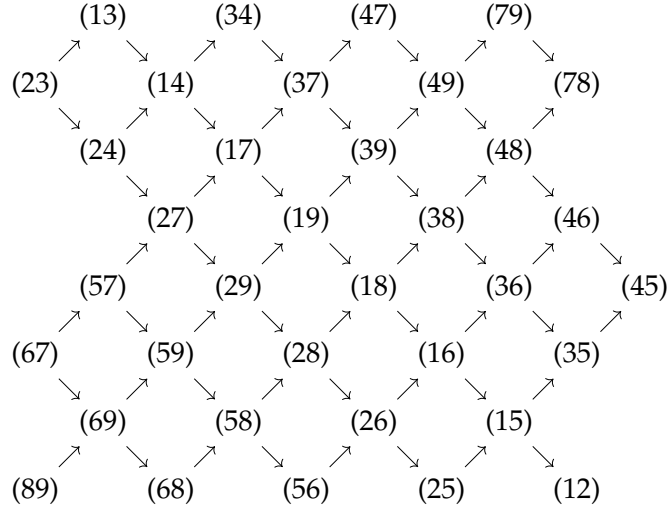

Let $c \in \mathfrak{S}_n$ be a Coxeter element. The \new{Auslander--Reiten quiver of $c$}, denoted $\AR(c)$, is the oriented graph satisfying the following conditions:
\begin{enumerate}[label = $\bullet$]
\item The vertices of $\AR(c)$ are the transpositions $(i,j)$, with $i<j$, in $\mathfrak{S}_n$;
\item The arrows of $\mathsf{AR}(c)$ are given, for all $i < j$, by
\begin{enumerate}[label = $\bullet$]
	\item $(i,j) \longrightarrow (i,c(j))$ whenever $i < c(j)$;
	
	\item $(i,j) \longrightarrow (c(i),j)$ whenever  $c(i) < j$.
\end{enumerate}
\end{enumerate}
To construct recursively such a graph, we can first find the initial adjacent transpositions of $c$, which are all the sources,  and step by step, using the second rule, construct the arrows and the vertices of $\AR(c)$ until we reach all the transpositions of $\mathfrak{S}_n$. Note that the sinks of $\AR(c)$ are given by the final adjacent transpositions of $c$.
See \cref{fig:ARc} for an explicit example.

\begin{rmk}
	The Auslander--Reiten quiver of a Coxeter element has a representation-theoretic meaning: briefly it corresponds to the oriented graph whose vertices are the indecomposable representations of a certain type $A$ quiver, and whose arrows are the irreducible morphisms between them. 
	
	To see further details about Auslander-Reiten quivers of type $A$ quivers in particular, we refer the reader to \cite[Section 3.1]{Sch14}. To learn more about quiver representation theory, and for more in-depth knowledge on the notion of Auslander--Reiten quivers, we invite the reader to look at \cite{ASS06}.
\end{rmk}

\section{An extended generalized Ferrers diagram RSK}

In the following, we describe a generalized version of RSK using (type $A$) Coxeter elements, and state the main result.

Let $\lambda$ be a nonzero integer partition and consider $(\B,\E)$ the unique elementary interval bipartition such that $\lambda(\B, \E) = \lambda$. Set $n = h_\lambda(1,1) + 1$. Let $c \in \mathfrak{S}_n$ and consider $\AR(c)$ its Auslander--Reiten quiver.

Recall that if $\B = \{b_1<\ldots<b_p\}$ and $\E = \{e_1 < \ldots < e_q \}$, then $(i,j) \in \Fer(\lambda)$ if and only if $b_i < e_{q-j+1}$. It allows us to label each box $(i,j)$ by a transposition $(b_i,e_{q-j+1}) \in \mathfrak{S}_n$. Thus it allows us to construct a one-to-one correspondence from fillings of shape $\lambda$ to fillings of the Auslander--Reiten quiver $\AR(c)$ which are supported on vertices $(b,e) \in \B \times \E$ such that $b < e$. Explicitly, for any filling $f$ of shape $\lambda$, we define $\overline{f}$ be the filling of $\AR(c)$ defined by $\overline{f}(b_i,e_{q-j+1}) = f(i,j)$ whenever $(i,j) \in \Fer(\lambda)$ and $\overline{f}(x,y) = 0$ otherwise.

As in \cref{s:Gans}, for $m \in \{1,\ldots, n-1\}$, let $(u_m,v_m)$ be the maximal pair with respect of $\unlhd$ in $D_m(\lambda)$. The boxes in the ideal generated by $(u_m, v_m)$ correspond to pairs $(i,j)$ such that $b_i \leqslant m < e_{q-j+1}$, and therefore $(u_m,v_m)$ is the maximal pair satisfying this condition.

For each $m \in \{1,\ldots,n-1\}$, we consider the subgraph $\AR_m(c)$ of $\AR(c)$ where the vertices are the transpositions $(i,j)$ with $i \leqslant m < j$. This subgraph has only one source and only one sink. 

We define $g = \RSK_{\lambda, c}(f)$ to be the fillings of shape $\lambda$ defined for $m \in \{1,\ldots,n-1\}$ by
\[\forall (i,j) \in D_m(\lambda),\quad g(i,j) = \GK_{\AR_m(c)}(f)_{u_m - i+1}.\]

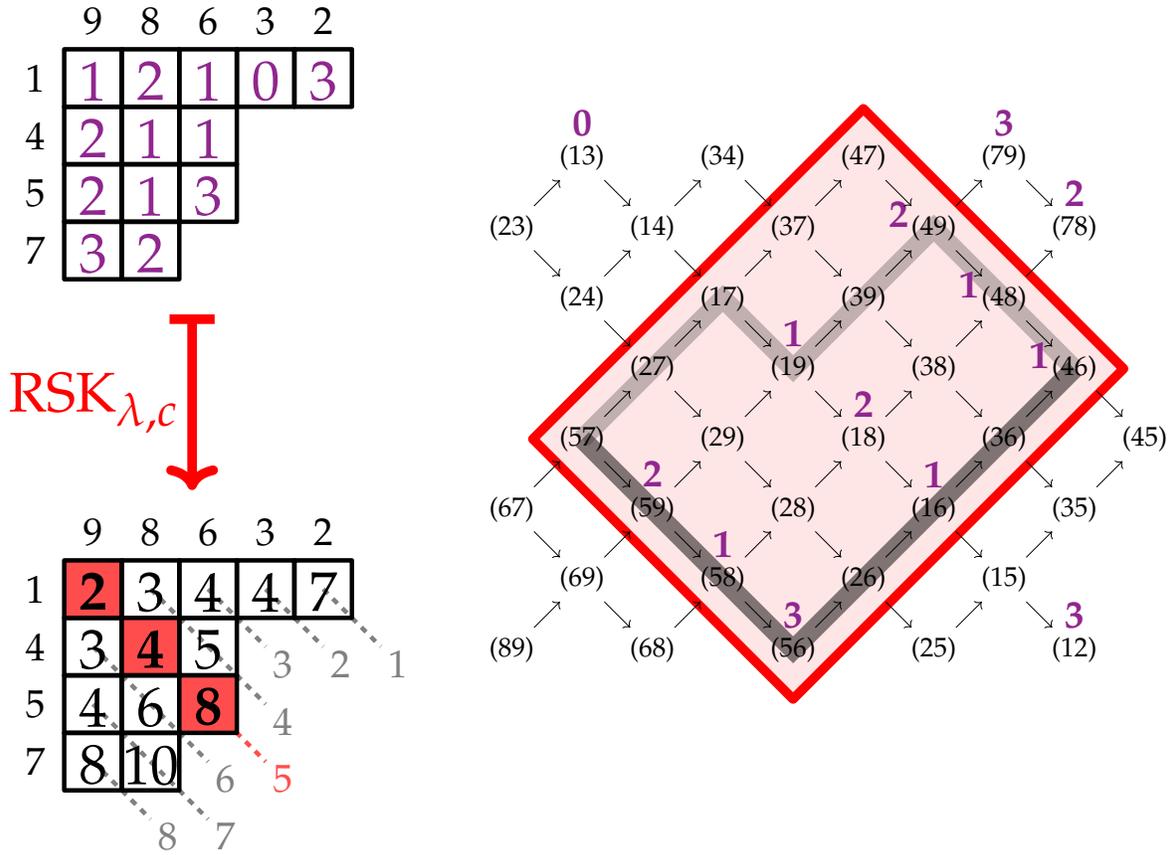
\begin{figure}[h!]
	\centering
	\[
	\scalebox{0.85}{
		\begin{tikzpicture}
			
			\begin{scope}[xshift=-7cm, yshift = 8.5cm,scale=0.9]
				\tkzDefPoint(0,0){a}
				\tkzDefPoint(0,1){b}
				\tkzDefPoint(1,1){c}
				\tkzDefPoint(1,0){d}
				\tkzDrawPolygon[line width = 0.7mm, color = black](a,b,c,d);
				
				\tkzDefPoint(1,0){a}
				\tkzDefPoint(1,1){b}
				\tkzDefPoint(2,1){c}
				\tkzDefPoint(2,0){d}
				\tkzDrawPolygon[line width = 0.7mm, color = black](a,b,c,d);
				
				\tkzDefPoint(2,0){a}
				\tkzDefPoint(2,1){b}
				\tkzDefPoint(3,1){c}
				\tkzDefPoint(3,0){d}
				\tkzDrawPolygon[line width = 0.7mm, color = black](a,b,c,d);
				
				\tkzDefPoint(3,0){a}
				\tkzDefPoint(3,1){b}
				\tkzDefPoint(4,1){c}
				\tkzDefPoint(4,0){d}
				\tkzDrawPolygon[line width = 0.7mm, color = black](a,b,c,d);
				
				\tkzDefPoint(4,1){a}
				\tkzDefPoint(4,0){b}
				\tkzDefPoint(5,0){c}
				\tkzDefPoint(5,1){d}
				\tkzDrawPolygon[line width = 0.7mm, color = black](a,b,c,d);
				
				\tkzDefPoint(0,0){a}
				\tkzDefPoint(0,-1){b}
				\tkzDefPoint(1,-1){c}
				\tkzDefPoint(1,0){d}
				\tkzDrawPolygon[line width = 0.7mm, color = black](a,b,c,d);
				
				\tkzDefPoint(1,0){a}
				\tkzDefPoint(1,-1){b}
				\tkzDefPoint(2,-1){c}
				\tkzDefPoint(2,0){d}
				\tkzDrawPolygon[line width = 0.7mm, color = black](a,b,c,d);
				
				\tkzDefPoint(2,0){a}
				\tkzDefPoint(2,-1){b}
				\tkzDefPoint(3,-1){c}
				\tkzDefPoint(3,0){d}
				\tkzDrawPolygon[line width = 0.7mm, color = black](a,b,c,d);
				
				\tkzDefPoint(0,-2){a}
				\tkzDefPoint(0,-1){b}
				\tkzDefPoint(1,-1){c}
				\tkzDefPoint(1,-2){d}
				\tkzDrawPolygon[line width = 0.7mm, color = black](a,b,c,d);
				
				\tkzDefPoint(1,-2){a}
				\tkzDefPoint(1,-1){b}
				\tkzDefPoint(2,-1){c}
				\tkzDefPoint(2,-2){d}
				\tkzDrawPolygon[line width = 0.7mm, color = black](a,b,c,d);
				
				\tkzDefPoint(2,-2){a}
				\tkzDefPoint(2,-1){b}
				\tkzDefPoint(3,-1){c}
				\tkzDefPoint(3,-2){d}
				\tkzDrawPolygon[line width = 0.7mm, color = black](a,b,c,d);
				
				\tkzDefPoint(0,-3){a}
				\tkzDefPoint(0,-2){b}
				\tkzDefPoint(1,-2){c}
				\tkzDefPoint(1,-3){d}
				\tkzDrawPolygon[line width = 0.7mm, color = black](a,b,c,d);
				
				\tkzDefPoint(1,-3){a}
				\tkzDefPoint(1,-2){b}
				\tkzDefPoint(2,-2){c}
				\tkzDefPoint(2,-3){d}
				\tkzDrawPolygon[line width = 0.7mm, color = black](a,b,c,d);
				
				\node at (0.5,1.5){{\Large $9$}};
				\node at (1.5,1.5){{\Large $8$}};
				\node at (2.5,1.5){{\Large $6$}};
				\node at (3.5,1.5){{\Large $3$}};
				\node at (4.5,1.5){{\Large $2$}};
				\node at (-0.5,0.5){{\Large $1$}};
				\node at (-0.5,-0.5){{\Large $4$}};
				\node at (-0.5,-1.5){{\Large $5$}};
				\node at (-0.5,-2.5){{\Large $7$}};
				
				\tkzLabelPoint[Plum](0.5,0.95){{\Huge $1$}};
				\tkzLabelPoint[Plum](1.5,0.95){{\Huge $2$}};
				\tkzLabelPoint[Plum](2.5,0.95){{\Huge $1$}};
				\tkzLabelPoint[Plum](3.5,0.95){{\Huge $0$}};
				\tkzLabelPoint[Plum](4.5,0.95){{\Huge $3$}};
				\tkzLabelPoint[Plum](0.5,-0.05){{\Huge $2$}};
				\tkzLabelPoint[Plum](1.5,-0.05){{\Huge $1$}};
				\tkzLabelPoint[Plum](2.5,-0.05){{\Huge $1$}};
				\tkzLabelPoint[Plum](0.5,-1.05){{\Huge $2$}};
				\tkzLabelPoint[Plum](1.5,-1.05){{\Huge $1$}};
				\tkzLabelPoint[Plum](2.5,-1.05){{\Huge $3$}};
				\tkzLabelPoint[Plum](0.5,-2.05){{\Huge $3$}};
				\tkzLabelPoint[Plum](1.5,-2.05){{\Huge $2$}};
			\end{scope}
			
			\begin{scope}[xshift=-7cm, yshift = 0.5cm,scale=0.9]
				
				\draw [line width=0.7mm, red!70, dashed] (0.5,0.5) --  (3.5,-2.5);
				\node[red!70] at (3.8,-2.8){{\Large $5$}};
				
				\draw [line width=0.7mm, gray, dashed] (1.5,0.5) --  (3.5,-1.5);
				\node[gray] at (3.8,-1.8){{\Large $4$}};
				
				\draw [line width=0.7mm, gray, dashed] (2.5,0.5) --  (3.5,-0.5);
				\node[gray] at (3.8,-0.8){{\Large $3$}};
				
				\draw [line width=0.7mm, gray, dashed] (3.5,0.5) --  (4.5,-0.5);
				\node[gray] at (4.8,-0.8){{\Large $2$}};
				
				\draw [line width=0.7mm, gray, dashed] (4.5,0.5) --  (5.5,-0.5);
				\node[gray] at (5.8,-0.8){{\Large $1$}};
				
				\draw [line width=0.7mm, gray, dashed] (0.5,-0.5) --  (2.5,-2.5);
				\node[gray] at (2.8,-2.8){{\Large $6$}};
				
				\draw [line width=0.7mm, gray, dashed] (0.5,-1.5) --  (2.5,-3.5);
				\node[gray] at (2.8,-3.8){{\Large $7$}};
				
				\draw [line width=0.7mm, gray, dashed] (0.5,-2.5) --  (1.5,-3.5);
				\node[gray] at (1.8,-3.8){{\Large $8$}};
				
				\tkzDefPoint(0,0){a}
				\tkzDefPoint(0,1){b}
				\tkzDefPoint(1,1){c}
				\tkzDefPoint(1,0){d}
				\tkzDrawPolygon[line width = 0.7mm, color = black, fill=red!70](a,b,c,d);
				
				\tkzDefPoint(1,0){a}
				\tkzDefPoint(1,1){b}
				\tkzDefPoint(2,1){c}
				\tkzDefPoint(2,0){d}
				\tkzDrawPolygon[line width = 0.7mm, color = black](a,b,c,d);
				
				\tkzDefPoint(2,0){a}
				\tkzDefPoint(2,1){b}
				\tkzDefPoint(3,1){c}
				\tkzDefPoint(3,0){d}
				\tkzDrawPolygon[line width = 0.7mm, color = black](a,b,c,d);
				
				\tkzDefPoint(3,0){a}
				\tkzDefPoint(3,1){b}
				\tkzDefPoint(4,1){c}
				\tkzDefPoint(4,0){d}
				\tkzDrawPolygon[line width = 0.7mm, color = black](a,b,c,d);
				
				\tkzDefPoint(4,1){a}
				\tkzDefPoint(4,0){b}
				\tkzDefPoint(5,0){c}
				\tkzDefPoint(5,1){d}
				\tkzDrawPolygon[line width = 0.7mm, color = black](a,b,c,d);
				
				\tkzDefPoint(0,0){a}
				\tkzDefPoint(0,-1){b}
				\tkzDefPoint(1,-1){c}
				\tkzDefPoint(1,0){d}
				\tkzDrawPolygon[line width = 0.7mm, color = black](a,b,c,d);
				
				\tkzDefPoint(1,0){a}
				\tkzDefPoint(1,-1){b}
				\tkzDefPoint(2,-1){c}
				\tkzDefPoint(2,0){d}
				\tkzDrawPolygon[line width = 0.7mm, color = black, fill = red!70](a,b,c,d);
				
				\tkzDefPoint(2,0){a}
				\tkzDefPoint(2,-1){b}
				\tkzDefPoint(3,-1){c}
				\tkzDefPoint(3,0){d}
				\tkzDrawPolygon[line width = 0.7mm, color = black](a,b,c,d);
				
				\tkzDefPoint(0,-2){a}
				\tkzDefPoint(0,-1){b}
				\tkzDefPoint(1,-1){c}
				\tkzDefPoint(1,-2){d}
				\tkzDrawPolygon[line width = 0.7mm, color = black](a,b,c,d);
				
				\tkzDefPoint(1,-2){a}
				\tkzDefPoint(1,-1){b}
				\tkzDefPoint(2,-1){c}
				\tkzDefPoint(2,-2){d}
				\tkzDrawPolygon[line width = 0.7mm, color = black](a,b,c,d);
				
				\tkzDefPoint(2,-2){a}
				\tkzDefPoint(2,-1){b}
				\tkzDefPoint(3,-1){c}
				\tkzDefPoint(3,-2){d}
				\tkzDrawPolygon[line width = 0.7mm, color = black,fill = red!70](a,b,c,d);
				
				\tkzDefPoint(0,-3){a}
				\tkzDefPoint(0,-2){b}
				\tkzDefPoint(1,-2){c}
				\tkzDefPoint(1,-3){d}
				\tkzDrawPolygon[line width = 0.7mm, color = black](a,b,c,d);
				
				\tkzDefPoint(1,-3){a}
				\tkzDefPoint(1,-2){b}
				\tkzDefPoint(2,-2){c}
				\tkzDefPoint(2,-3){d}
				\tkzDrawPolygon[line width = 0.7mm, color = black](a,b,c,d);
				
				\node at (0.5,1.5){{\Large $9$}};
				\node at (1.5,1.5){{\Large $8$}};
				\node at (2.5,1.5){{\Large $6$}};
				\node at (3.5,1.5){{\Large $3$}};
				\node at (4.5,1.5){{\Large $2$}};
				\node at (-0.5,0.5){{\Large $1$}};
				\node at (-0.5,-0.5){{\Large $4$}};
				\node at (-0.5,-1.5){{\Large $5$}};
				\node at (-0.5,-2.5){{\Large $7$}};
				
				\tkzLabelPoint(0.5,0.95){{\Huge $\mathbf{2}$}};
				\tkzLabelPoint(1.5,0.95){{\Huge $3$}};
				\tkzLabelPoint(2.5,0.95){{\Huge $4$}};
				\tkzLabelPoint(3.5,0.95){{\Huge $4$}};
				\tkzLabelPoint(4.5,0.95){{\Huge $7$}};
				\tkzLabelPoint(0.5,-0.05){{\Huge $3$}};
				\tkzLabelPoint(1.5,-0.05){{\Huge $\mathbf{4}$}};
				\tkzLabelPoint(2.5,-0.05){{\Huge $5$}};
				\tkzLabelPoint(0.5,-1.05){{\Huge $4$}};
				\tkzLabelPoint(1.5,-1.05){{\Huge $6$}};
				\tkzLabelPoint(2.5,-1.05){{\Huge $\mathbf{8}$}};
				\tkzLabelPoint(0.5,-2.05){{\Huge $8$}};
				\tkzLabelPoint(1.5,-2.05){{\Huge $10$}};
			\end{scope}
			
			\begin{scope}[yshift=0cm,scale=1.1]
			\tkzDefPoint(0.3,3){a}
			\tkzDefPoint(4,-0.7){b}
			\tkzDefPoint(8.7,4){c}
			\tkzDefPoint(5,7.7){d}
			\tkzDrawPolygon[line width = 1.5mm, color = red, fill = red!10](a,b,c,d);
			
			\draw[line width = 3mm,  opacity=0.5, black] (1,3) -- (4,0) -- (8,4);
			
			\draw[line width = 3mm, opacity=0.3, black!80] (1,3) -- (3,5) -- (4,4) -- (6,6) -- (8,4);

			\node (67) at (0,2) {(67)};
			\node (57) at (1,3) {(57)};
			\node (27) at (2,4) {(27)};
			\node (17) at (3,5) {(17)};
			\node (37) at (4,6) {(37)};
			\node (47) at (5,7) {(47)};
			
			\node (69) at (1,1) {(69)};
			\node (59) at (2,2) {(59)};
			\node[Plum] at (2,2.5){{\Large $\mathbf 2$}};
			\node (29) at (3,3) {(29)};
			\node (19) at (4,4) {(19)};
			\node[Plum] at (4,4.5){{\Large $\mathbf 1$}};
			
			\node (39) at (5,5) {(39)};
			
			\node (49) at (6,6) {(49)};
			\node[Plum] at (5.5,6.2){{\Large $\mathbf 2$}};
			
			\node (68) at (2,0) {(68)};
			\node (58) at (3,1) {(58)};
			\node[Plum] at (3,1.5){{\Large $\mathbf 1$}};
			\node (28) at (4,2) {(28)};
			\node (18) at (5,3) {(18)};
			\node[Plum] at (5,3.5){{\Large $\mathbf 2$}};
			\node (38) at (6,4) {(38)};
			\node (48) at (7,5) {(48)};
			\node[Plum] at (6.5,5.2){{\Large $\mathbf 1$}};
			\node (78) at (8,6) {(78)};
			\node[Plum] at (8,6.5){{\Large $\mathbf 2$}};
			\node (89) at (0,0) {(89)};
			\node (79) at (7,7) {(79)};
			\node[Plum] at (7,7.5){{\Large $\mathbf 3$}};
			
			\node (56) at (4,0) {(56)};
			\node[Plum] at (4,0.5){{\Large $\mathbf 3$}};
			
			\node (25) at (6,0) {(25)};
			\node (12) at (8,0) {(12)};
			\node[Plum] at (8,0.5){{\Large $\mathbf 3$}};
			
			\node (13) at (1,7) {(13)};
			\node[Plum] at (1,7.5){{\Large $\mathbf 0$}};
			\node (34) at (3,7) {(34)};
			
			\node (26) at (5,1) {(26)};
			\node (15) at (7,1) {(15)};
			\node (23) at (0,6) {(23)};
			\node (14) at (2,6) {(14)};
			
			\node (16) at (6,2) {(16)};
			\node[Plum] at (6,2.5){{\Large $\mathbf 1$}};
			
			\node (35) at (8,2) {(35)};
			\node (24) at (1,5) {(24)};
			
			\node (36) at (7,3) {(36)};
			\node (45) at (9,3) {(45)};
			
			\node (46) at (8,4) {(46)};
			\node[Plum] at (7.5,4.2){{\Large $\mathbf 1$}};
			
			\draw[->] (13)--(14);\draw[->] (14)--(34);\draw[->] (14)--(17);\draw[->] (15)--(35);\draw[->] (15)--(12);\draw[->] (16)--(36);\draw[->] (16)--(15);\draw[->] (17)--(37);\draw[->] (17)--(19);\draw[->] (18)--(38);\draw[->] (18)--(16);\draw[->] (19)--(39);\draw[->] (19)--(18);\draw[->] (23)--(13);\draw[->] (23)--(24);\draw[->] (24)--(14);\draw[->] (24)--(27);\draw[->] (25)--(15);\draw[->] (26)--(16);\draw[->] (26)--(25);\draw[->] (27)--(17);\draw[->] (27)--(29);\draw[->] (28)--(18);\draw[->] (28)--(26);\draw[->] (29)--(19);\draw[->] (29)--(28);\draw[->] (34)--(37);\draw[->] (35)--(45);\draw[->] (36)--(46);\draw[->] (36)--(35);\draw[->] (37)--(47);\draw[->] (37)--(39);\draw[->] (38)--(48);\draw[->] (38)--(36);\draw[->] (39)--(49);\draw[->] (39)--(38);\draw[->] (46)--(45);\draw[->] (47)--(49);\draw[->] (48)--(78);\draw[->] (48)--(46);\draw[->] (49)--(79);\draw[->] (49)--(48);\draw[->] (56)--(26);\draw[->] (57)--(27);\draw[->] (57)--(59);\draw[->] (58)--(28);\draw[->] (58)--(56);\draw[->] (59)--(29);\draw[->] (59)--(58);\draw[->] (67)--(57);\draw[->] (67)--(69);\draw[->] (68)--(58);\draw[->] (69)--(59);\draw[->] (69)--(68);\draw[->] (79)--(78);\draw[->] (89)--(69);
			\end{scope}	
			\draw [|->,line width=1.5mm,red] (-5,5.25) -- node[left]{{\Huge $\RSK_{\lambda,c} $}} (-5,2.5);	
	\end{tikzpicture} }\]
\caption{\label{fig:RSKandAR} Explicit calculation of $\RSK_{\lambda,c}(f)$ for the boxes in $D_5(\lambda)$ from a filling of $\lambda = (5,3,3,2)$, with $c = (1,3,4,7,9,8,6,5,2)$}
\end{figure}

See \cref{fig:RSKandAR} for an explicit example.

Our main result is the following.

\begin{thm} \label{mainthm}
	Let $\lambda$ be a nonzero integer partition. Consider $n = h_\lambda(1,1) + 1$. Let $c \in \mathfrak{S}_n$ be a Coxeter element. The map $\RSK_{\lambda,c}$ gives a one-to-one correspondence from fillings of shape $\lambda$ to reverse plane partitions of shape $\lambda$.
\end{thm}

The following result shows that we extended the RSK correspondence.

\begin{prop}
	Let $\lambda$ be a nonzero integer partition. Consider $n = h_\lambda(1,1) + 1$ and $(\B,\E)$ be the only elementary interval bipartition such that $\lambda(\B,\E) = \lambda$. Let $c \in \mathfrak{S}_n$ be the Coxeter element such that
	\begin{enumerate}[label = $\bullet$]
		\item for $i \in \{1,\ldots, n-1\}$, $(i,i+1)$ is final in $c$ if and only if $i \in \B$ and $i+1 \in \E$;
		
		\item for $i \in \{2, \ldots, n-2\}$, $(i,i+1)$ is initial in $c$ if and only if $i \in \E$ and $i+1 \in \B$. 
	\end{enumerate}
	Then $\RSK_{\lambda,c} = \RSK_\lambda$. Moreover, $c$ and $c^{-1}$ are the unique Coxeter element of $\mathfrak{S}_n$ satisfying this property.
\end{prop}

\begin{rmk}
	Gansner's RSK for a fixed integer partition $\lambda$ admits a local description in terms of toggles on $G_\lambda$. Based on the proof given in \cite{Deq23}, for $c=(1,2,\ldots,n)$, we can give a local description in terms of toggles on $\AR(c)$. However, more works need to be done for a general choice of $c$, as this local description does not extend naturally.  
\end{rmk}

\section{Some words about quiver representation theory}
\label{s:quiver}

This section aims to give a dictionary to link the result from \cite{Deq23} with \cref{mainthm}. 

Fix $Q=(Q_0,Q_1)$ a type $A$ quiver. A \new{finite dimensional representation $E$ of $Q$ over $\mathbb{C}$} is an assignment of a finite dimensional $\mathbb{C}$-vector space $E_q$ to each vertex $q$ of $Q$, and an assignment of a $\mathbb{C}$-linear transformation $E_\alpha : E_i \longrightarrow E_j$ to each arrow $\alpha : i \rightarrow j$ of $Q$. For two representations $E$ and $F$, a morphism $\phi : E \longrightarrow F$ is the data of a $\mathbb{C}$-linear map $\phi_q$ for each vertex $q$ of $Q$ such that for any arrow $\alpha : i \rightarrow j$, $\phi_j E_\alpha = F_\alpha \phi_i$.  Denote by $\rep_\mathbb{K}(Q)$ the category consisting of the representations of $Q$. 

Any representation $E$ of $Q$ can be uniquely decomposed into a direct sum of copies of indecomposable representations up to isomorphism. Thus, we can consider the invariant which counts the number of indecomposable summands of each isomorphism class in $E$. Write it $\Mult(E)$.

In \cite{GPT19}, A. Garver, R. Patrias and H. Thomas introduce a new invariant of quiver representations, called the generic Jordan form data. For any representation $E$ of $Q$, write $\GenJF(E)$ for the generic Jordan form data of $E$. This data encodes the generic behavior of a nilpotent endomorphism $N = (N_q)_{q \in Q_0}$ of the representation via the size of the Jordan blocks of each $N_q$. In some subcategories, the representation can be recovered from this invariant up to isomorphism.

They also show that the map from $\Mult$ to $\GenJF$ generalizes the RSK correspondence for type A quivers,
using Gansner’s previous work \cite{Ga81Ac}. 

As this map is bijective, if we restrict it to the representation in some subcategories $\mathscr{C}$, one can be interested to get an explicit way to invert it. An algebraic method developed in \cite{GPT19} asks the subcategory $\mathscr{C}$ to satisfy the following property. For any $E \in \mathscr{C}$,  there exists a dense open set $\Omega$ (in the Zariski topology) in the set of representations admitting a nilpotent endomorphism with Jordan forms encoded by $\GenJF(E)$ such that any $F \in \Omega$ is isomorphic to $E$. Such a subcategory is said to be \new{canonically Jordan recoverable (CJR)}.

More recently, in \cite{Deq23}, we gave a combinatorial characterization of all the CJR subcategories of representations of $Q$, substancially enlarging the family of subcategories for which $\GenJF$ is a complete invariant given by \cite{GPT19}. The maximal such subcategories can be described thanks to the elementary interval partitions $(\B,\E)$ of $\{1, \ldots, n+1\}$.

The following table compares the representation-theoretic tools used in \cite{Deq23} and the combinatorial tools used to describe our generalized RSK.

\begin{center}
	\begin{tabular}{|c|c|}
	\hline \textbf{Combinatorial tools} & \textbf{Representation-theoretic tools} \\
	
	\hline Coxeter element of $\mathfrak{S}_n$ & Orientation of an $A_{n-1}$ type quiver $Q$ \\
	
	\hline Transposition in $\mathfrak{S}_n$ & Indecomposable representation in $\rep_\mathbb{C}(Q)$ \\
	
	\hline AR quiver of $c$ & AR quiver of $\rep_\mathbb{C}(Q)$ \\
	
	\hline Integer partition $\lambda$ with $h_\lambda(1,1) = n-1$ & CJR subcategory $\mathscr{C}$ of $\rep_\mathbb{C}(Q)$ \\
	
	\hline Filling of $\lambda$ &  $\Mult(E)$ for some $E \in \mathscr{C}$ \\
	
	\hline Reverse plane partition of $\lambda$ & $\GenJF(E)$ for $E$ in $\mathscr{C}$. \\
	\hline 
\end{tabular}
\end{center}

\acknowledgements{I acknowledge the Institut des Sciences Mathématiques of Canada for its partial funding support. I thank my Ph.D. supervisor, Hugh Thomas, for all our discussions on this subject, his helpful advice, and his support throughout my thesis work, and even more.}


\printbibliography

\end{document}